\documentclass[onecolumn]{autart3}    

\usepackage{graphicx}      
\usepackage{amssymb,amsmath}
\usepackage{dsfont}
\usepackage{color}
\usepackage{epstopdf}        
\usepackage{amsfonts}
\usepackage{enumerate}
\usepackage{picins}

\theoremstyle{definition}
\newtheorem{definition}{Definition}[section]

\newtheorem{remark}{Remark}[section]
\newtheorem{proposition}{Proposition}[section]

\newtheorem{lemma}{Lemma}[section]

\newtheorem{assumption}{Assumption}[section]
\newcommand*{\QEDB}{\hfill\ensuremath{\square}}

\usepackage{color}

\newcommand{\tcbb}{\textcolor{black}}
\begin{document}
\begin{frontmatter}
\title{Robust Hybrid Zero-Order Optimization Algorithms with Acceleration via Averaging in Time\thanksref{footnoteinfo}}
%
\thanks[footnoteinfo]{Corresponding author: J. I. Poveda. \tcbb{This work was supported in part by the grants NSF CRII - CNS - 1947613, NSF CAREER 1553407, AFOSR YIP: FA9550-18-1-0150, and ONR YIP: N00014-19-1-2217.}} 
\author[Colorado]{Jorge I. Poveda}\ead{jorge.poveda@colorado.edu},    
\author[Harvard]{Na Li}\ead{nali@seas.harvard.edu}               

\address[Colorado]{Department of Electrical, Computer, and Energy Engineering, University of Colorado, Boulder, CO, 80309 USA.}  
\address[Harvard]{School of Engineering and Applied Sciences, Harvard University, Cambridge, MA, 02138 USA.}     


\maketitle

\begin{abstract}
This paper presents a new class of robust zero-order algorithms for the solution of real-time optimization problems with acceleration. In particular, we propose a family of extremum seeking (ES) dynamics that can be universally modeled as singularly perturbed hybrid dynamical systems with restarting mechanisms.  From this family of dynamics, we synthesize four fast algorithms for the solution of convex, strongly convex, constrained, and unconstrained optimization problems. In each case, we establish robust semi-global practical asymptotic or exponential stability results, and we also show how to obtain well-posed discretized algorithms that retain the main properties of the original dynamics. Given that existing averaging theorems for singularly perturbed hybrid systems are not directly applicable to our setting, we derive an extended averaging theorem that relaxes some of the assumptions made in the literature, allowing us to make a clear link between the $\mathcal{K}\mathcal{L}$ bounds that characterize the rates of convergence of the hybrid dynamics and their average dynamics. We also show that our results are applicable to non-hybrid algorithms, thus providing a general framework for accelerated ES dynamics based on averaging theory. We present different numerical examples to illustrate our results. 
\end{abstract}
\end{frontmatter}
%
%
%
%
\section{Introduction}
This paper studies feedback-based algorithms for the solution of optimization problems of the form
\begin{equation}\label{problem}
\min~~\phi(z)~~~\text{subject to}~~~z\in \mathcal{F},
\end{equation}
where $\phi:\mathbb{R}^n\to\mathbb{R}$ is a smooth cost function and $\mathcal{F}\subset\mathbb{R}^n$ is a nonempty, closed, and convex set. Unlike most of the standard model-based optimization problems considered in the literature, our main assumption is that the mathematical forms of the cost function and its gradient are unknown, and the algorithms have access only to real-time evaluations of the cost. Feedback-based optimization algorithms of this form, also called extremum seeking (ES) dynamics \cite{KrsticBookESC,tan06Auto,TAC13NesicESC,Guay:03}, have recently seen a renewed interest, and several novel architectures and theoretical results have been presented in \cite{DerivativesESC,ScheinkerBoundedLetters,Grushkovskaya2017,Suttner2017} and \cite{ThiagoKrstic} for ordinary differential equations (ODEs), and in \cite{Poveda:16,BlackBoxHybridES,Ronny_Kutadinata} for systems modeled as hybrid dynamical systems. However, while significant progress has been made during the last years, one of the persistent challenges in ES is how to guarantee fast rates of convergence without sacrificing stability and robustness properties that are critical for real-world implementations in noisy environments. In the model-based optimization domain, this challenge has driven the development of several \emph{accelerated} optimization algorithms that incorporate dynamic momentum, see for instance \cite{InexactNesterov2,ODE_Nesterov,Wibisono1,Mihailo19,Faziyab18Siam,HamiltonianHybrid}, and references therein. When the gradient of the cost function is not available, accelerated discrete-time algorithms have also been studied in \cite{InexactNesterov2} and \cite{Nesterov_derivative_free}. In the continuous-time domain, ES dynamics based on adaptations of the \emph{time-invariant} Heavy-ball method were studied in \cite{HeavyBollES} and \cite{Michalowsky15}. However, in the context of ES, no acceleration properties have been established so far, and the tradeoffs that may emerge between robustness and acceleration remain mostly unexplored.

Motivated by this background, in this paper we introduce the first \emph{accelerated} ES algorithms with dynamic momentum based on averaging theory. These ES algorithms can be seen as model-free versions of continuous-time Nesterov's gradient dynamics with discrete-time restarting mechanisms, and therefore they are naturally modeled as hybrid dynamical systems (HDS). As recently shown in \cite{PovedaNaliCDC19,HamiltonianHybrid,Candes_Restarting}, restarting mechanisms are instrumental in order to regularize and/or improve the stability and transient properties of accelerated time-varying optimization dynamics that otherwise may not be suitable for applications in feedback-based optimization. Indeed, unlike existing results in the literature of ES based on the \emph{time-invariant} Heavy-Ball method \cite{HeavyBollES,Michalowsky15}, the restarting mechanisms used in this paper allow us to exploit the underlaying acceleration properties of the  \emph{time-varying} Nesterov's ordinary differential equations (ODEs) studied in \cite{ODE_Nesterov,Wibisono1,HighResolution2018}, without sacrificing the robustness and stability properties that are critical in ES. This feature further allows us to establish semi-global practical asymptotic stability results with $\mathcal{K}\mathcal{L}$ bounds for all our algorithms, as well as fast rates of convergence for convex and strongly convex functions in constrained and unconstrained optimization problems. Moreover, since averaging theory for HDS can also be applied to ODEs, our results are also applicable to continuous-time ES dynamics. We exploit this property in order to establish additional novel (semi-global practical) \emph{exponential} stability results for Primal-Dual and Augmented Primal Dual ES algorithms.  The dynamics considered in this paper are modeled and analyzed using the framework of set-valued HDS presented in \cite{Goebel:12}. By using this framework, we establish novel structural robustness results for all our algorithms, as well as convergence results for discretized dynamics obtained via Euler or Runge-Kutta discretization, which extend previous discretization results \cite{EnenbauerZeroOrder} to more general hybrid settings. Since existing averaging theorems for set-valued HDS require an average system with a  \emph{uniformly globally asymptotically stable} (UGAS) compact set, a condition that is generally not satisfied in ES, we extend the averaging theorem to systems having average dynamics with \emph{semi-global practical asymptotic stability} properties. This auxiliary result is a modest extension of the results of \cite{averagingTeel}, \cite{Wang:12_Automatica}, \cite{Poveda:16}, that allows us to directly link the $\mathcal{K}\mathcal{L}$ function of the average HDS with the $\mathcal{K}\mathcal{L}$ function of the original dynamics. Similar relaxations have been considered in \cite{tan06Auto} for ODEs, and in \cite[Thm. 1]{Ronny_Kutadinata} for HDS with \emph{non-hybrid} average systems. 

The rest of this paper is organized as follows: Sections \ref{sec_preliminaries} and \ref{sec_algorithms} present the preliminaries and main results. Sections \ref{structural_robustness} and \ref{section_numeric} present discretization results and numerical examples. Section \ref{averaging_section} presents the averaging framework needed for the analysis of the algorithms. Section \ref{sec_analysis} presents all the proofs, and finally Section \ref{sec_conclusions} ends with some conclusions.
%
\section{Preliminaries}
\label{sec_preliminaries}
The set of (nonnegative) real numbers is denoted by ($\mathbb{R}_{\geq0}$) $\mathbb{R}$. We use $\mathbb{B}$ to denote a closed unit ball of appropriate dimension, $\rho\mathbb{B}$ to denote a closed ball of radius $\rho>0$, and $\mathcal{X}+\rho\mathbb{B}$ to denote the union of all sets obtained by taking a closed ball of radius $\rho$ around each point in the set $\mathcal{X}$. The closed convex hull of a set $\mathcal{X}$ is denoted as $\overline{\text{co}}(\mathcal{X})$, and we use $\lim \sup_{i\to\infty}\mathcal{X}_i$ to denote the outer limit of a sequence of sets $\{\mathcal{X}_i\}^{\infty}_i$ \cite[Def. 5.1]{Goebel:12}. We use $\mathbb{S}^1\subset\mathbb{R}^2$ to denote the unit circle centered at the origin, and $\mathbb{T}^n:=\mathbb{S}^1\times\mathbb{S}^1\times\ldots\times\mathbb{S}^1$ to denote the $n$-Cartesian product of $\mathbb{S}^1$. Given a vector $x\in\mathbb{R}^n$ and a compact set $\mathcal{A}$ we use $|x|_{\mathcal{A}}:=\min_{y\in\mathcal{A}}\|y-x\|_2$. We also use $I_n\in\mathbb{R}^{n\times n}$ to denote the identity matrix, $\mathbf{c}_n\in\mathbb{R}^n$ to denote the vector with all entries equal to $c\in\mathbb{R}$, and $e_i$ to denote a unit vector of appropriate dimension with $i^{th}$ entry equal to $1$.  A continuous function $\alpha:\mathbb{R}_{\geq0}\to\mathbb{R}_{\geq0}$ is of class $\mathcal{K}$ if $\alpha$ is zero at zero and strictly increasing. It is said to be of class $\mathcal{K}_{\infty}$ if it is of class $\mathcal{K}$ and grows unbounded. A function $\sigma:\mathbb{R}_{\geq0}\to\mathbb{R}_{\geq0}$ is of class $\mathcal{L}$ if it is continuous, non-increasing, and converging to zero as its argument grows unbounded. A function $\beta$ is of class $\mathcal{K}\mathcal{L}$ if it is of class $\mathcal{K}$ in its first argument, and of class $\mathcal{L}$ in its second argument. A function $\phi:\mathbb{R}^n\to\mathbb{R}$ is said to be radially unbounded if $\phi(x)\to\infty$ as $|x|\to\infty$, and it is said to be of class $\mathcal{C}^k$ if its $k^{th}$ derivative is continuous. 

To study our optimization algorithms, we use the framework of HDS \cite{Goebel:12}, which considers systems of the form
\begin{subequations}\label{HDS}
\begin{align}
&x\in C,~~~~~~~~\dot{x}= F(x),~~~\label{flowseq1}\\
&x\in D,~~~~~~x^+= G(x),~~~\label{jumpseq1}
\end{align}
\end{subequations}
where $x\in\mathbb{R}^n$ is the state, $F:\mathbb{R}^n\to\mathbb{R}^n$ is called the flow map, and $G:\mathbb{R}^n\to\mathbb{R}^n$ is called the jump map. The sets $C$ and $D$, called the flow set and the jump set, respectively, characterize the points in the space where the system evolves according to \eqref{flowseq1}, or \eqref{jumpseq1}, respectively. The \emph{data} of the HDS is defined as $\mathcal{H}:=\{C, F, D, G\}$. Systems of the form \eqref{HDS} generalize continuous-time systems and discrete-time systems. Namely, continuous-time systems can be seen as HDS of the form \eqref{HDS} with $D=\emptyset$, while discrete-time systems correspond to the case when $C=\emptyset$. Solutions $x:\text{dom}(x)\to\mathbb{R}^n$ to \eqref{HDS} are defined on hybrid time domains, and they are parametrized by a continuous-time index $t\in\mathbb{R}_{\geq0}$ and a discrete-time index $j\in\mathbb{Z}_{\geq0}$. Solutions with an unbounded time domain are said to be complete. For a precise definition of solutions to HDS we refer the reader to Appendix B.
\begin{definition}\label{definition1}
A HDS  \eqref{HDS} is said to satisfy the Basic Conditions  if $C$ and $D$ are closed, $C\subset\text{dom}(F)$, $D\subset\text{dom}(G)$, and $F$ and $G$ are continuous on $C$ and $D$, respectively. 
\end{definition}
To study systems of the form \eqref{HDS}, we will use the following definitions that are standard in hybrid systems \cite{Goebel:12}.
\begin{definition}\label{definitionUGAS}
A compact set $\mathcal{A}\subset\mathbb{R}^n$ is said to be Uniformly Globally pre-Asymptotically Stable (UGpAS) for  a HDS $\mathcal{H}$ if there exists $\beta\in\mathcal{K}\mathcal{L}$ such that every solution $x$ of $\mathcal{H}$ satisfies $|x(t,j)|_{\mathcal{A}}\leq \beta(|x(0,0)|_{\mathcal{A}},t+j)$, for all $(t,j)~\in\text{dom}(x)$. When $\beta(r,s)=c_1r\exp(-c_2s)$ for some $c_1,c_2>0$, we say that $\mathcal{H}$ renders $\mathcal{A}$ Uniformly Globally pre-Exponentially Stable (UGpES). If, additionally, all solutions are complete, we use the acronyms UGAS and UGES, respectively.
\end{definition}
\begin{definition}\label{definitionSGPAS}
For a parameterized HDS $\mathcal{H}_{\delta_1,\delta_2}$ a compact set $\mathcal{A}\subset\mathbb{R}^n$ is said to be Semi-Globally Practically pre-Asymptotically Stable (SGPpAS) as $(\delta_1,\delta_2)\to0^+$ with $\beta\in\mathcal{K}\mathcal{L}$ if for all compact sets $K\subset\mathbb{R}^n$ and all $\nu>0$,  $\exists$ $\delta_1^*>0$ such that $\forall$ $\delta_1\in(0,\delta_1^*)$,  $\exists$ $\delta_2^*>0$ such that $\forall$ $\delta_2\in(0,\delta_2^*)$, every solution $x_{\delta_1,\delta_2}$ of $\mathcal{H}_{\delta_1,\delta_2}$ with $x_{\delta_1,\delta_2}(0,0)\in K$ satisfies
\begin{equation}\label{wdelta}
|x_{\delta_1,\delta_2}(t,j)|_{\mathcal{A}}\leq \beta(|x_{\delta_1,\delta_2}(0,0)|_{\mathcal{A}},t+j)+\nu,
\end{equation}
\vspace{0.1cm}
for all $(t,j)~\in\text{dom}(x_{\delta_1,\delta_2})$.  When $\beta$ has exponential form we say that $\mathcal{A}$ is Semi-Globally Practically pre-Exponentially Stable (SGPpES) as $(\delta_1,\delta_2)\to0^+$ with $\beta\in\mathcal{K}\mathcal{L}$. If, additionally, all solutions satisfying the bound \eqref{wdelta} are complete, we use the acronyms SGPAS and SGPES, respectively.
\end{definition}
When $D=\emptyset$, and all solutions are complete, Definition \ref{definitionUGAS} reduces to the standard UGAS and UGES notions for continuous-time systems. Also, when $\mathcal{H}_{\delta}$ is parameterized by only one constant $\delta>0$, Definition \ref{definitionSGPAS} recovers the standard definitions of SGPAS and SGPES.
\begin{remark}
In Definition \ref{definitionSGPAS} the order of the parameters $\delta_1$ and $\delta_2$ is relevant, i.e., in general $\delta_2^*$ may depend on $\delta_1$ in a non-trivial way. Definition \ref{definitionSGPAS} can be extended to parameterized HDS with any number of parameters, i.e., $\mathcal{H}_{\delta_1,\delta_2,\ldots,\delta_p}$, $p\in\mathbb{Z}_{\geq1}$. 
\end{remark}
In this paper, we are interested in optimization algorithms with desirable \emph{robustness} properties with respect to small disturbances that are unavoidable in practice. The following definition aims to capture this property.
\begin{definition}\label{def_robustness}
Let $\mathcal{H}$ render UGpAS (resp. SGPpAS as $\delta\to0^+$) a compact set $\mathcal{A}$ with $\beta\in\mathcal{K}\mathcal{L}$.  We say that $\mathcal{H}$ is Structurally Robust if for all measurable functions $e:\mathbb{R}_{\geq0}\to\mathbb{R}^n$ satisfying $\sup_{t\geq0}|e(t)|\leq \bar{e}$, with $\bar{e}>0$, the perturbed system
\begin{subequations}\label{perturbed_system}
\begin{align}
&x+e\in C,~~~~~~~~~\dot{x}=F(x+e)+e,\\
&x+e\in D,~~~~~~~x^+= G(x+e)+e,
\end{align}
\end{subequations}
renders the set $\mathcal{A}$ SGPpAS as $\bar{e}\to0^+$ (resp. SGPpAS as $(\delta,\bar{e})\to 0^+$) with $\beta\in \mathcal{K}\mathcal{L}$.
\end{definition}
As noted in \cite{PovedaNaliCDC19}, some continuos-time accelerated gradient dynamics may not satisfy the robustness property of Definition \ref{def_robustness}. This tradeoff between robustness and acceleration makes it not trivial to design robust accelerated ES algorithms by using standard averaging tools for ODEs.
%
\section{Hybrid Accelerated Extremum Seeking: Algorithms and Main Stability Results}
\label{sec_algorithms}
Let the set of solutions of problem \eqref{problem} be given by
\begin{equation}\label{optimal_set}
\mathcal{A}_{\phi}:=\left\{z^*\in\mathcal{F}:\phi(z^*)\leq \phi(z),~\forall~z\in\mathcal{F}\right\},
\end{equation}
where $\phi^*:=\phi(\mathcal{A}_{\phi})>-\infty$. We consider a family of hybrid accelerated extremum seeking (HAES) algorithms that can be modeled as a HDS with states $x\in\mathbb{R}^{n+m},~\tau\in\mathbb{R}_{>0}$ and $\mu\in\mathbb{R}^{2n}$, with the following data
\begin{equation}\label{HAES}
\mathcal{H}_{a,\varepsilon}=\{C_{es},F_{es},D_{es},G_{es}\}, 
\end{equation}
where $\tau$ models a restarting timer, and $\mu$ models an excitation signal. The state $x:=[x_1^\top,x_2^\top]\in\mathbb{R}^{n+m}$ has two main components, with $x_1\in\mathbb{R}^n$ acting as the main state, and $x_2\in\mathbb{R}^m$ acting as an auxiliary state that is instrumental for the incorporation of dynamic momentum or dual variables in the algorithms. The continuous-time dynamics of system \eqref{HAES} are parameterized by two tunable positive constants: $a\in\mathbb{R}_{>0}$ and $\varepsilon\in\mathbb{R}_{>0}$. These continuous-time dynamics are characterized by the following flow set and flow map
%
\begin{subequations}\label{flows_main_equations}
\begin{align}
(x,\tau,\mu)&\in C_{es}:= \mathbb{R}^{n+m}\times \mathcal{T}_C \times \mathbb{T}^n,\label{flow_main}\\
\left(\begin{array}{c}
\dot{x}\\
\dot{\tau}\\
\dot{\mu}
\end{array}\right)&=F_{es}(x,\tau,\mu):=
\left(\begin{array}{c}
F_x(x,\mu,\phi(z))\\
F_{\tau} \\
\dfrac{1}{\varepsilon}R\mu
\end{array}\right),\label{flows_generalHDS}
\end{align}
\end{subequations}
where the mappings $(F_x, F_{\tau})$ and the set $\mathcal{T}_C\subset\mathbb{R}_{>0}$ will be designed based on the qualitative assumptions made on the cost function $\phi$. The dynamics of the state $\mu$ are characterized by the matrix $R\in\mathbb{R}^{2n\times2n}$, which is a block diagonal matrix with ${\ell}^{th}$ diagonal block given by $R_{\ell}:=2\pi\kappa_\ell\cdot\left[-e_2,e_1\right]\in\mathbb{R}^{2\times2}$, where $\kappa_{\ell}\in\mathbb{R}_{>0}$ is a tunable parameter. These linear dynamics describe $n$ uncoupled oscillators that generate solutions $\mu:\mathbb{R}_{\geq0}\to\mathbb{R}^{2n}$ with odd entries
\begin{equation}\label{sinusoid2}
\mu_{i}(t)=\Psi_{i}(t)^\top \mu_{0,i}~~~~i\in\{1,3,5,\ldots,2n-1\},
\end{equation}
where $\Psi_{i}(t):=\left[\cos\left(\frac{2\pi t}{\varepsilon}\kappa_{\frac{i+1}{2}}\right),\sin\left(\frac{2\pi t}{\varepsilon}\kappa_{\frac{i+1}{2}}\right)\right]^\top$ and $\mu_{0,i}:=[\mu_{i}(0),\mu_{i+1}(0)]^\top$. The argument $z$ of the cost function $\phi$ in equation \eqref{flows_generalHDS} is updated via the feedback law
\begin{equation}\label{input_1}
z=x_1+a\tilde{\mu},~~~~\tilde{\mu}:=[\mu_1,\mu_3,\mu_5,\ldots,\mu_{2n-1}]^\top,
\end{equation}
where $\mu_i$ is given by \eqref{sinusoid2}. We will make the following assumption on the parameters $\kappa_i$ to guarantee suitable averaging properties for the signal $\mu$.
\begin{assumption}\label{assumption_freq}
For each $\ell\in\{1,2,\ldots,n\}$ the parameter $\kappa_{\ell}$ is a positive rational number, and $\kappa_{\ell}\neq \kappa_j$ for all $j\neq \ell$. 
\end{assumption}
The discrete-time dynamics of system \eqref{HAES} are characterized by the following jump set and jump map:
\begin{subequations}\label{main_jump_dynamics}
\begin{align}
(x,\tau,\mu)&\in D_{es}:=\mathbb{R}^{n+m}\times \mathcal{T}_D \times \mathbb{T}^n.\label{jump_main}\\
\left(\begin{array}{c}
x^+\\
\tau^+\\
\mu^+
\end{array}\right)&=G_{es}(x,\tau,\mu):=
\left(\begin{array}{c}
G_x(x)\\
T_{\min}\\
\mu
\end{array}\right),~\label{jumps_main}
\end{align}
\end{subequations}
where $T_{\min}\in\mathbb{R}_{>0}$. The mapping $G_x$ and the set $\mathcal{T}_D\subset\mathbb{R}_{>0}$ will also be designed based on the qualitative assumptions made on the cost function $\phi$.  


In order to study in a unified manner the stability and convergence properties of system $\mathcal{H}_{a,\varepsilon}$, we define the set
\begin{equation}\label{optimal}
\mathcal{A}:=\mathcal{A}_x\times\mathcal{T}_C\times\mathbb{T}^n,
\end{equation}
where $\mathcal{A}_x\subset\mathbb{R}^{n+m}$ is a closed and bounded set having the property that its projection onto $\mathbb{R}^n$ coincides with the set of solutions of \eqref{problem}, i.e., 
\begin{equation}\label{projection}
\left\{x_1\in\mathbb{R}^n: x=[x_1^\top,x_2^\top]^\top\in \mathcal{A}_x\right\}=\mathcal{A}_{\phi}.
\end{equation} 
Based on these definitions, our goal is to design the mappings $(F_x,F_\tau, G_x)$ and the sets $(\mathcal{T}_C,\mathcal{T}_D)$ to guarantee suitable stability and fast convergence properties with respect to the set $\mathcal{A}$ for the HDS \eqref{HAES}. In order to do this, we will focus on four main qualitative optimization problems of the form \eqref{problem}. Namely: 1) convex cost functions with no constraints; 2) strongly convex functions with no constraints; 3) strongly convex functions with equality constraints; and 4) strongly convex functions with inequality constraints. 
\subsection{Case 1: Unconstrained Convex Optimization}
\label{subsection:1}
We first consider the case when $\mathcal{F}:=\mathbb{R}^n$, and the cost function $\phi$ satisfies the following assumption:
\begin{assumption}\label{assumption1}
The mapping $z\mapsto \phi(z)$ is $\mathcal{C}^2$, convex, radially unbounded, and satisfies at least one of the following conditions: (a) $\phi$ has a unique minimizer; (b) $\nabla \phi$ is globally Lipschitz.
\end{assumption}
For functions $\phi$ satisfying Assumption \ref{assumption1}, we consider a HAES $\mathcal{H}_{a,\varepsilon}$ with a constant restarting frequency $F_{\tau}>0$, and mappings:
\begin{equation}\label{maps1}
F_x:=\left(\begin{array}{c}
\dfrac{2}{\tau}(x_2-x_1)-\dfrac{2}{a}k_1\phi(z)\tilde{\mu}\\
-\dfrac{4}{a}k_2\tau \phi(z)\tilde{\mu}
\end{array}\right),~~~G_x:=\left(\begin{array}{c}x_1\\x_2\end{array}\right),
\end{equation}
with sets
\begin{equation}\label{sets_1}
\mathcal{T}_C:=[T_{\min},~T_{\max}],~~~~\mathcal{T}_D:=[T_{\text{med}},~~T_{\max}],
\end{equation}
where $m=n$, $a\in\mathbb{R}_{>0}$ is the same parameter of \eqref{input_1}, $k_1\in\mathbb{R}_{\geq0}$ and $k_2\in \mathbb{R}_{>0}$ are tunable gains, $z$ is given by the feedback law \eqref{input_1}, and $T_{\max},T_{\text{med}}\in\mathbb{R}_{>0}$ satisfy $T_{\max}\geq T_{\text{med}}$ and $T_{\text{med}}-T_{\min}>\epsilon>0$, for some $\epsilon\in\mathbb{R}_{>0}$. This hybrid system incorporates dynamic momentum during the flows via the state $x_2$, and it can generate non-unique solutions from a given initial condition, including solutions with periodic and aperiodic restarting. In particular, the HAES allows jumps whenever $\tau\geq T_{\text{med}}$ but no later than when $\tau=T_{\max}$. Since $\tau^+=T_{\min}\notin \mathcal{T}_D$ and $\epsilon>0$, any two consecutive jumps in a given solution of the system are separated at least by a positive amount of time $\epsilon/F_{\tau}$ during which the system has to flow. Thus, every solution is uniformly non-Zeno. For the case when $T_{\max}=T_{\text{med}}$ the jumps (i.e., restartings) are periodic. Additionally, since $G_{es}(D_{es})\subset C_{es}\cup D_{es}$, by item (c) in Lemma \ref{existence_completness} in the Appendix, the solutions of the HDS do not stop due to jumps or flows leaving the set $C_{es}\cup D_{es}$. Indeed, by construction, the HAES satisfies the Basic Conditions of Definition \ref{definition1}.  

The following theorem characterizes the stability, acceleration, and robustness properties of the  HAES   with respect to the compact set \eqref{optimal}, with $\mathcal{A}_x$ defined as
\begin{equation}\label{optimal_1}
\mathcal{A}_x:=\{x\in\mathbb{R}^{2n}:x_1=x_2,~x_1\in\mathcal{A}_{\phi}\}.
\end{equation}
\tcbb{Below, we express the convergence bound on the sub-optimality measure in terms of a Lyapunov function $V_{k_1}$ defined in Section \ref{sec_analysis} for the \emph{average} hybrid dynamics of system \eqref{HAES}.} 
\begin{thm}\label{thm1a}
Suppose that Assumptions \ref{assumption_freq} and \ref{assumption1} hold, and consider the HAES \eqref{HAES} with state $\tilde{x}:=[x^\top,\mu^\top,\tau]^\top$ and data given by \eqref{maps1} and \eqref{sets_1}. Then, the following holds with $k_1=0$ and $F_{\tau}=\frac{1}{2}$:
\begin{enumerate}[(a)]
\item The set $\mathcal{A}$ is SGPAS as $(a,\varepsilon)\to0^+$ with $\beta_1\in\mathcal{K}\mathcal{L}$. Additionally, system $\mathcal{H}_{a,\varepsilon}$ is Structurally Robust. 
 \item \tcbb{For each compact set $K_0\subset\mathbb{R}^{2n}$ such that $\mathcal{A}_x\subset \text{int}(K_0)$, and each $\nu>0$,  $\exists$  $a^*>0$ such that $\forall~a\in(0,a^*)$, $\exists$ $\varepsilon^*>0$ such that $\forall~\varepsilon\in(0,\varepsilon^*)$, all solutions with $x(0,0)\in K_0$ induce the bound:}
\begin{equation}\label{inequality_algo1}
\tcbb{\phi(z(t,j))-\phi^*\leq \frac{4V_{k_1}(\tilde{x}(\underline{t}_j,j))}{k_2(t-\underbar{t}_{j})^2}+\nu,}
\end{equation}
for all $(t,j)\in \text{dom}(\tilde{x})$ \tcbb{such that $t>\underbar{t}_{j}$, where $\underbar{t}_{j}=\min\{t: (t,j)\in\text{dom}(\tilde{x})\}$, and $V_{k_1}(\cdot)$ is a Lyapunov function for the average hybrid system of $\mathcal{H}_{a,\varepsilon}$ that satisfies $\lim\sup_{j\to\infty} V_{k_1}(\tilde{x}(\underline{t}_j,j))\leq \nu$.}
\end{enumerate}
Moreover, if the minimizer of $\phi$ is unique, items (a)-(b) also hold with $k_1\geq0$ and $F_{\tau}=1$. 
\end{thm}
In words, item (a) of Theorem \ref{thm1a} establishes that by \emph{orderly} tuning the parameters $a$ and $\varepsilon$, system $\mathcal{H}_{a,\varepsilon}$ guarantees robust convergence of the state $x$, and therefore $z$ via \eqref{input_1}, to any arbitrarily small neighborhood of the set of minimizers $\mathcal{A}_{\phi}$. Note that the robustness margins $\bar{e}$ may in general depend on the parameters $a$ and $\varepsilon$. \tcbb{On the other hand, item (b) describes a semi-acceleration property during the flows that is novel in the literature of ES. Namely, for each $j\in\mathbb{Z}_{\geq0}$ and all $t>\underline{t}_{j}$ such that $(t,j)\in\text{dom}(\tilde{x})$, the sub-optimality measure $\phi(z)-\phi^*$ will decrease at a rate of $\mathcal{O}_j(1/(t-\underline{t}_j)^2)$, modulo a small residual error. Since for $j=0$ we have that $\underline{t}_j=0$, by using $\tau(0,0)=T_{\text{min}}$  the right-hand side of \eqref{inequality_algo1} simplifies to $c_0/t^2+\nu$ for all $t>0$ in the interval $(0,F^{-1}_{\tau}(T_{\text{med}}-T_{\min})]$, where $c_0=(4/k_2)V_{k_1}(\tilde{x}(0,0))$ is a constant defined by the initial conditions of the algorithm.} Given that this interval of flow can be made arbitrarily large by the choice of $T_{\text{med}}$, for ``flat'' convex cost functions the acceleration property \eqref{inequality_algo1} can induce an initial dramatical improvement in the rate of convergence of the ES dynamics in comparison to standard gradient descent-based ES algorithms, which, in general, only achieve rates of convergence of order $\mathcal{O}(1/t)$ for the class of smooth convex cost functions \cite[pp. 7]{ODE_Nesterov}. 
\subsubsection*{Tuning Guidelines and Connections with Nesterov's ODE}
\label{remark1}
In order to achieve the semi-acceleration property \eqref{inequality_algo1}, the mapping $F_x$ in \eqref{maps1} is designed to be intrinsically related to the time-varying accelerated Nesterov's ODE \cite{ODE_Nesterov,HighResolution2018}. Indeed, as shown in the analysis of Section \ref{section_average_systems}, by using the change of variables $s=x_1$ and $x_2=s+0.5\tau\left(\dot{s}+k_1\nabla\phi(x_1)\right)$, the $x$-component of the solutions of $\mathcal{H}_{a,\varepsilon}$ approximates the behavior of its \emph{average system}, which can be written as
\begin{equation}\label{NesterovsODE}
\ddot{s}+\frac{(2+\dot{\tau})\dot{s}}{\tau}+4k_2\nabla \phi(s)+k_1\left(\nabla^2\phi(s)^\top\dot{s}+\frac{\dot{\tau}}{\tau}\nabla\phi(s)\right)=0.
\end{equation}
When $k_1=0$, $k_2=1$, $\dot{\tau}=0.5$ and $\tau(0)\geq1$, equation \eqref{NesterovsODE} corresponds to the ODE studied in \cite{Jadbabaie_RK}. When $k_1=0$ and $\dot{\tau}=0$, equation \eqref{NesterovsODE} reduces to the time-invariant Heavy-Ball dynamics, studied in the context of ES in \cite{HeavyBollES}. When $k_1=1/\sqrt{L}$, $\dot{\tau}=1$ and $k_2=0.25$ equation \eqref{NesterovsODE} corresponds to the Hessian-driven dynamics studied in \cite{AcceleratedMethods19}. Interestingly, in this case the average dynamics \eqref{NesterovsODE} incorporate the Hessian matrix of $\phi$ even though there is no explicit Hessian estimation in the extremum seeking dynamics \eqref{maps1}. Finally, note that when $k_1=0$, $\dot{\tau}=1$, $k_2=0.25$, $T_{\min}=0$, $T_{\text{med}}=\infty$, and $\dot{x}(0)=\tau(0)=0$, equation \eqref{NesterovsODE} reduces to the time-varying Nesterov's ODE studied in \cite{ODE_Nesterov}.  As shown in \cite[Ex. 1]{PovedaNaliCDC19}, for this dynamics the persistent restarting of $\tau$ is needed to induce uniform convergence, which is closely related to structural robustness. For this reason, there are clear tradeoffs between robustness and acceleration in the hybrid system $\mathcal{H}_{a,\varepsilon}$. Namely, as $T_{\text{med}}\to\infty$ the intervals of flows satisfying \eqref{inequality_algo1} grow larger (i.e., less frequent restarting), but the robustness margins $\bar{e}$ of the perturbed system \eqref{perturbed_system} shrink to zero. In the limit, when $T_{\text{med}}=\infty$, system $\mathcal{H}_{a,\varepsilon}$ behaves as the time-varying Nesterov's ODE, which is highly sensitive to arbitrarily small disturbances. Thus, the tuning of the parameter  $T_{\text{med}}$ is critical in order to obtain a good tradeoff between longer periods of flow with acceleration \eqref{inequality_algo1} and larger margins of robustness. Finally, when $k_1>0$ the HAES incorporates a Hessian-driven damping term that has been shown to slightly improve the transient performance of Nesterov's ODE \cite{HighResolution2018,AcceleratedMethods19}. However, good performance can be obtained in the hybrid dynamics $\mathcal{H}_{a,\varepsilon}$ even when $k_1=0$. Illustrative numerical examples are presented in Section \ref{section_numeric}.
%
\subsection{Case 2: Unconstrained Strongly Convex Optimization}
\label{subsection:2}
We now study accelerated ES dynamics with \emph{momentum restarting} mechanisms that induce (semi-global practical) exponential stability with rates of convergence adjustable by the restarting frequency. 

We consider cost functions $\phi$ that satisfy the following assumption:
\begin{assumption}\label{assumption2}
The mapping $z\mapsto\phi(z)$ is $\mathcal{C}^2$, and there exist $\theta>0$ and $L>0$ such that $|\nabla \phi(z')-\nabla \phi(z'')|\leq L|z'-z''|$ and $(\nabla \phi(z')-\nabla\phi(z''))^\top (z'-z'')\geq \theta|z'-z''|^2$, for all $z',z''\in\mathbb{R}^n$.
%
%
%
\end{assumption}
Under Assumption \ref{assumption2}, the set of minimizers of $\phi$ is a singleton, i.e., $\mathcal{A}_{\phi}=\{z^*\}$. In this case, we consider HAES dynamics with the following mappings:
\begin{equation}\label{maps2}
F_x:=\left(\begin{array}{c}
\dfrac{2}{\tau}(x_2-x_1)\\
-\dfrac{4}{a}k\tau \phi(z)\tilde{\mu}
\end{array}\right),~~F_{\tau}=\frac{1}{2},~~G_x:=\left(\begin{array}{c}
x_1\\
x_1
\end{array}\right),
\end{equation}
where $m=n$, $k\in\mathbb{R}_{>0}$ is a tunable gain, and 
\begin{equation}\label{sets_2}
\mathcal{T}_C:=[T_{\min},~~T_{\max}],~~~~~\mathcal{T}_D:=\{T_{\max}\},
\end{equation}
with $T_{\max}-T_{\min}>\epsilon$ for some $\epsilon>0$. By construction this HDS also satisfies the Basic Conditions, and its jumps are periodic and separated by an interval of length $\Delta T:=2(T_{\max}-T_{\min})> 2\varepsilon>0$. Thus, every solution is uniformly non Zeno. The following theorem characterizes the stability, acceleration, and robustness properties of the HAES $\mathcal{H}_{a,\varepsilon}$ with respect to the set $\mathcal{A}$ in \eqref{optimal} with $\mathcal{A}_x=\{z^*\}\times\{z^*\}$. 
\begin{thm}\label{thm2}
Suppose that Assumptions \ref{assumption_freq} and \ref{assumption2} hold, and consider the HAES \eqref{HAES} with state $\tilde{x}:=[x^\top,\mu^\top,\tau]^\top$ and data given by \eqref{maps2} and \eqref{sets_2}.  Let the parameters ($k,T_{\min},T_{\max}$) satisfy the inequality:
\begin{equation}\label{dwell_time}
T^2_{\max}-T_{\min}^2\geq \frac{1}{2\theta k}.
\end{equation}
Then, the following holds:
\begin{enumerate}[(a)]
\item The compact set $\mathcal{A}$ is SGPES as $(a,\varepsilon)\to0^+$ with $\beta_2\in\mathcal{K}\mathcal{L}$. Additionally, system $\mathcal{H}_{a,\varepsilon}$ is Structurally Robust. 
\item \tcbb{For each compact set $K_0\subset\mathbb{R}^{2n}$ such that $\mathcal{A}_x\subset \text{int}(K_0)$, and each $\nu>0$,  $\exists$  $a^*>0$ such that $\forall~a\in(0,a^*)$, $\exists$ $\varepsilon^*>0$ such that $\forall~\varepsilon\in(0,\varepsilon^*)$, all solutions with $x(0,0)\in K_0$, $\tau(0,0)=T_{\min}$, and $x_1(0,0)=x_2(0,0)$, satisfy the bound}
\begin{equation}\label{exponential_decrease_during_jumps}
\tcbb{\phi(z(t,j))-\phi^*\leq \alpha_0 \tilde{\gamma}^j\big(\phi(z(0,0))-\phi^*\big)+\nu,}
\end{equation}
\tcbb{for all $(t,j)\in\text{dom}(\tilde{x})$ such that $t>\underline{t}_j:=\min\{t:(t,j)\in\text{dom}(\tilde{x})\}$, where  $\tilde{\gamma}:=\frac{1}{kT_{\max}^2}\left(\frac{1}{2\theta}+kT_{\min}^2\right)$ and $\alpha_0:=T_{\max}^2/T_{\min}^2$.}
\end{enumerate}
\end{thm}
Item (a) of Theorem \ref{thm2} says that condition \eqref{dwell_time} is sufficient for semi-global practical exponential convergence of the states $x_1$ and $x_2$ to the point $z^*$.  When $T_{\min}+T_{\max}>1$, condition \eqref{dwell_time} can be satisfied by a standard \emph{dwell-time} condition of the form $T_{\max}-T_{\min}>(2k\theta)^{-1}$, which has not been studied before in the stability analysis of ES controllers. \tcbb{On the other hand, since condition \eqref{dwell_time} implies that $\tilde{\gamma}$ in \eqref{exponential_decrease_during_jumps} satisfies $\tilde{\gamma}\in(0,1)$, item (b) establishes an explicit constant decrease of the sub-optimality measure during jumps, modulo a small residual $\nu$-error. This property, induced by the restarting mechanism, is novel in the context of averaging-based ES algorithms, and for certain classes of cost functions it can be further exploited to achieve faster rates of convergence compared to gradient descent-based ES methods.}
\subsubsection*{Tuning Guidelines: Quasi-Optimal Restarting vs. Black-Box Restarting}
\label{tuningstrong}
\tcbb{For cost functions $\phi$ satisfying Assumption \ref{assumption2}, and for a fixed gain $k>0$, and any $\nu>0$, the gradient descent-based ES algorithms of the form $\dot{x}_1=-k\frac{2}{a}\phi(z)\tilde{\mu}$, with $z=x_1+a\tilde{\mu}$, generate convergence bounds on the sub-optimality measure of the form $\phi(z(t))-\phi^*\leq (\phi(z(0))-\phi^*)e^{-2k\theta t}+\frac{\nu}{2}$,  provided $(a,\varepsilon)$ are selected sufficiently small (cf. Theorem 2) \cite{KrsticBookESC,DerivativesESC}. 
Therefore, a $\nu$-error in the sub-optimality measure is achieved when $t\geq t_G^*:=(1/2k\theta )\log((\phi(z(0))-\phi^*)/0.5\nu)$. On the other hand, for the same constants $k>0$ and $\nu>0$, and by using an appropriate choice of $T_{\max}$ and small values of $T_{\min}$, the HAES can exploit the restarting dynamics to induce convergence times that are approximately of order $\mathcal{O}(1/\sqrt{k\theta}\log(1/\nu))$. For example, if knowledge of $\theta$ is available, this can be achieved by using}
\begin{equation}\label{quasi_optimal}
\tcbb{T^*_{\max}=e\sqrt{\frac{1}{2k\theta}+T_{\min}^2},}
\end{equation}
\tcbb{which leads to $\tilde{\gamma}=1/e^2$ in \eqref{exponential_decrease_during_jumps}. In this case, the inequality $\alpha_0\tilde{\gamma}^j\big(\phi(z(0,0))-\phi^*\big)\leq 0.5\nu$ holds whenever $j\geq 0.5 \log(\alpha_0(\phi(z(0,0))-\phi^*)/0.5\nu)$. Multiplying both sides of the inequality by the switching period $\Delta T=2(T^*_{\max}-T_{\min})$, we obtain that $\phi(t,j)-\phi^*\leq \nu$ for all times $t\geq t^*_\mathcal{H}$ with}
\begin{equation*}
\tcbb{t^*_\mathcal{H}:=\left(e\sqrt{\frac{1}{2k\theta}+T_{\min}^2}-T_{\min}\right)\log\left(\frac{\alpha_0(\phi(0,0)-\phi^*)}{0.5\nu}\right).}
\end{equation*}
\tcbb{For small (but fixed) values of $T_{\min}>0$ we obtain that $t^*_\mathcal{H}$ is approximately of order $\mathcal{O}(1/\sqrt{k\theta}\log(\frac{1}{\nu}))$. Note that when $T^2_{\min}\approx 0$, the expression \eqref{quasi_optimal} can be explicitly computed by using $\Delta T=2T_{\max}$ and by minimizing $\tilde{\gamma}^{c/\Delta T}$ over $T_{\max}$, for any $c>0$, which leads to the ``optimal'' restarting period $\Delta T^*=e\sqrt{2/k\theta}$. Similar optimal periodic restarting conditions have been established in the discrete-time optimization literature \cite[Sec. 3]{Candes_Restarting}, and in some accelerated continuous-time algorithms  \cite{PovedaNaliCDC19}, \cite{ODE_Nesterov}. However, to the knowledge of the authors this type of result has not been established before in the context of ES. If, additionally, we set $k=1/2L$ in both the HAES and the gradient descent-based ES dynamics, we recover the well-known convergence bounds of order $\mathcal{O}(\sqrt{L/\theta}\log(1/\nu))$ and $\mathcal{O}(L/\theta~\log(1/\nu))$, respectively.  Since the HAES requires $T_{\min}>0$, we refer to \eqref{quasi_optimal} as a ``quasi-optimal'' restarting condition.}
\begin{remark}
\tcbb{While the existence of $T^*_{\max}$ highlights a theoretical advantage of the HAES over gradient-descent ES for strongly convex cost functions with $\theta\ll1$ (or large condition numbers if $k=1/2L$), it is important to note that in ES problems the values of the constants $\theta$ and $L$ are usually unknown, and therefore it is difficult to use in practice the exact
restarting parameter \eqref{quasi_optimal}. However, numerical examples show that similar rates of convergence can be obtained by using ``black-box'' restarting parameters $(T_{\min},T_{\max})$ obtained after a few tuning iterations. This black-box property has also been observed in the discrete-time optimization literature, e.g., \cite{Candes_Restarting,restarting_black_box}.} The restarting mechanism of the HAES gives an extra degree of freedom to tune the convergence properties of the algorithm. In Section \ref{section_numeric} we present illustrative numerical examples for different values of $T_{\max}$. 
\end{remark}

%
\subsection{Case 3: Strongly Convex Optimization with Linear Equality Constraints}
\label{subsection:3}
The framework of zero-order HDS can also be used to study algorithms that do not have jumps, i.e., dynamics modeled as ODEs. In this section, we illustrate this idea by considering an optimization problem \eqref{problem} with a feasible set given by
\begin{equation}\label{feasible_saddle1}
\mathcal{F}:=\{z\in \mathbb{R}^n:Az=b\},~~b\in\mathbb{R}^{m},~~A\in\mathbb{R}^{m\times n},
\end{equation}
which satisfies the following assumption:
\begin{assumption}\label{assumption3}
The matrix $A$ is full row rank, and there exist $\upsilon_1,\upsilon_2>0$ such that  $\upsilon_1I_m\leq AA^\top\leq\upsilon_2I_m$. 
\end{assumption}
For Case 3, we consider a HAES with the the following mappings:
\begin{equation}\label{map3}
F_x:=\left(\begin{array}{c}
-\dfrac{2}{a}\phi(z)\tilde{\mu}-kA^\top x_2  \\
(Ax_1-b)
\end{array}\right),~~~G_x:=x,
\end{equation}
where $k\in\mathbb{R}_{>0}$ is a tunable gain. Since these dynamics are independent of $\tau$, we can define the mapping $F_{\tau}$ and the sets $\mathcal{T}_C$ and $\mathcal{T}_D$ as follows: 
\begin{equation}\label{sets_3}
F_{\tau}=0,~~~~\mathcal{T}_C:=\{T_{\min}\},~~~~~~\mathcal{T}_D:=\{\emptyset\},
\end{equation}
which generates an empty jump set $D_{es}$ in \eqref{jump_main}, i.e., system $\mathcal{H}_{a,\varepsilon}$ does not generate solutions that jump. Nevertheless, the resulting HAES still satisfies the Basic Conditions and can be studied using the same tools as in Cases 1 and 2. To analyze the stability and convergence properties of the algorithm, we consider the set
\begin{align}\label{optimal_lag1}
\mathcal{A}_x:=&\Big\{(x^*_1,x^*_2)\in\mathbb{R}^{n+m}:\mathcal{L}(x^*_1,x_2)\leq\mathcal{L}(x^*_1,x^*_2)\leq\mathcal{L}(x_1,x^*_2),\notag\\
&~~~~~~~\forall~x_1\neq x^*_1,~x_2\neq x_2^*\Big\},
\end{align}
which is the set of saddle points of the Lagrangian 
\begin{equation}\label{lagrangian1}
\mathcal{L}(x_1,x_2)=\phi(x_1)+x_2^\top (Ax_1-b).
\end{equation}
Under Assumptions \ref{assumption2} and \ref{assumption3}, and by the results of \cite{GuananLetters}, the set $\mathcal{A}_x$ is a singleton and satisfies \eqref{projection}. 
\begin{thm}\label{thm3}
Consider the set $\mathcal{A}$ given by \eqref{optimal} with $\mathcal{A}_x$ given by \eqref{optimal_lag1} with Lagrangian \eqref{lagrangian1}. Suppose that Assumptions \ref{assumption_freq}, \ref{assumption2}, and \ref{assumption3} hold, and consider the HAES  \eqref{HAES} with data given by \eqref{map3} and \eqref{sets_3}.  Then, the compact set $\mathcal{A}$ is SGPES as $(a,\varepsilon)\to0^+$ with $\beta_3\in\mathcal{K}\mathcal{L}$. Additionally, system $\mathcal{H}_{a,\varepsilon}$ is Structurally Robust. 
%
%
%
%
\end{thm}
%
%
The dynamics characterized by equation \eqref{map3} can be seen as a type of zero-order Primal-Dual algorithm \cite{GuananLetters}, similar to those considered in \cite{Durr13L,Ye:16b}, for which \emph{asymptotic} convergence results have been established. However, Theorem \ref{thm3} shows that the semi-global practical convergence result is indeed \emph{exponential} and Structurally Robust.

%
%
%
\subsection{Case 4: Strongly Convex Cost Function with Inequality Constraints}
\label{subsection:4}
We finish this section by considering strongly convex functions with inequality constraints:
\begin{equation}
\mathcal{F}:=\{z\in\mathbb{R}^n: Az\leq b\},~~A\in\mathbb{R}^{m\times n},~~b\in\mathbb{R}^{m}.
\end{equation}
We consider the mappings $F_x$ and $G_x$ given by
\begin{equation}\label{map4}
F_x:=\left(\begin{array}{c}
-\dfrac{2}{a}\phi(z)\mu-k\sum_{j=1}^mH_j(x)A_j\\
\sum_{j=1}^{\tcbb{m}}\left(H_j(x)-x_{2,j}\right)e_j
\end{array}\right),~~G_x:=\{x\},
\end{equation}
where $A^\top:=[A_1,A_2,\ldots,A_m]$, $A_i\in\mathbb{R}^n$ for all $i$, $b:=[b_1,b_2,\ldots,b_m]^\top$, $H_j(x):=\max\left(A^\top_jx_1-b_j+x_{2,j},0\right)$, $k\in\mathbb{R}_{>0}$ is a tunable gain, and the data ($F_{\tau}$,$\mathcal{T}_C$, $\mathcal{T}_D$) are again given by \eqref{sets_3}, i.e., the solutions of the system do not jump. In this case, the mapping $F_x$ in \eqref{map4} describes a class of novel \emph{augmented Primal-Dual extremum seeking dynamics}, defined with respect to the following augmented Lagrangian:
\begin{equation}\label{lagrangian2}
\mathcal{L}(x_1,x_2):=\phi(x_1)+\sum_{j=1}^m H_j(x),
\end{equation}
which, under strong convexity of $\phi$ and Assumption \ref{assumption3}, also generates a singleton set $\mathcal{A}_x$ given by \eqref{optimal_lag1}, \cite{GuananLetters}. The following theorem also establishes a robust semi-global practical exponential stability result for system $\mathcal{H}_{a,\varepsilon}$.
\begin{thm}\label{thm4}
Consider the set $\mathcal{A}$ given by \eqref{optimal} generated from the set $\mathcal{A}_x$ given by \eqref{optimal_lag1} with Augmented  Lagrangian \eqref{lagrangian2}. Suppose that Assumptions \ref{assumption_freq}, \ref{assumption2}, and \ref{assumption3} hold and consider the  HAES  \eqref{HAES} with data given by \eqref{sets_3} and \eqref{map4}. Then, the compact set $\mathcal{A}$ is SGPES as $(a,\varepsilon)\to0^+$ with $\beta_4\in\mathcal{K}\mathcal{L}$. Moreover, system $\mathcal{H}_{a,\varepsilon}$ is Structurally Robust.  
%
%
\end{thm}
%
%
%
%
\section{Stable Discretization via Forward-Euler and Consistent Runge-Kutta}
\label{structural_robustness}
\label{sect_discretization}
Hybrid systems satisfying the Basic Conditions and having suitable asymptotic stability properties are ``stable'' under Euler and Runge-Kutta discretization. Motivated by this fact, as well as by recent discretization results in accelerated optimization \cite{Jadbabaie_RK}, we now study how to construct a suitable discretization for the HAES \eqref{HAES} with a fixed step size $h>0$. The resulting zero-order discretized system, which has overall state $\bar{x}_h=(x_h,\tau_h,\mu_h)$, and discrete-time dynamics
\begin{subequations}\label{discretized_HANDS}
\begin{align}
&\bar{x}_h\in C_h,~~~~~\bar{x}_h^+=F_h(\bar{x}_h),\label{discretize_flows}\\
&\bar{x}_h\in D_h,~~~~~\bar{x}_h^+=G_h(\bar{x}_h),
\end{align}
\end{subequations}
will retain the $\mathcal{K}\mathcal{L}$ convergence bounds of the original zero-order hybrid dynamics \eqref{HAES}, up to a time scaling of $t=\ell h$, where $\ell$ is the discrete-time index of the discretized flows \eqref{discretize_flows}. In order to do this, we rely on the notion of well-posed hybrid simulators, introduced in \cite{SimulatorHybridSystems}.
\begin{definition}\label{regular_discretization}
The discretized HAES obtained from $\mathcal{H}_{a,\varepsilon}$ and denoted by $\mathcal{H}_h$, is said to be \emph{well-posed} if the discretized data $(F_h,C_h,G_h,D_h)$ satisfies the following conditions:
\begin{enumerate}[(a)]
\item  $F_h$ is such that, for each compact set $K\subset\mathbb{R}^{n+m+1+2n}$, there exists a function $\rho\in\mathcal{K}_\infty$ and $h^*>0$ such that for each $\bar{x}_h\in C_h\cap K$ and each $h\in(0,h^*]$
\begin{equation*}
F_h(\bar{x}_h)\in \bar{x}_h+h~\overline{\text{co}} \left(F_{es}\big(\bar{x}_h+\rho(h)\mathbb{B}\big)\right)+h\rho(h)\mathbb{B}.
\end{equation*}
\item $G_h$ is such that for any decreasing sequence $h_i\to0^+$ we have that $G_0=G_{es}$, where $G_0$ is the graphical limit \cite[Def. 5.18]{Goebel:12} of $G_{h_i}$ as $h_i\to0^+$.
\item The sets $C_h$ and $D_h$ are such that for any sequence $\{h_i\}^{\infty}_{i=1}\searrow0$ such that $h_i\to 0^+$ we have that $\lim\sup_{i\to\infty}~C_{h_i}\subset C_{es}$ and $\lim\sup_{i\to\infty}~D_{h_i}\subset D_{es}$, \cite[Def. 5.1]{Goebel:12}.  
\end{enumerate}
\end{definition}
Let $\sum_{k=1}^{\bar{s}} b_k=1$, and $S=\{1,2,\ldots,\bar{s}\}$, $\bar{s}\in\mathbb{Z}_{>1}$. Then, by \cite[Ex. 4.8 $\&$ 4.9]{SimulatorHybridSystems} the forward-Euler method $F_h(\bar{x}_h):=\bar{x}_h+hF(\bar{x}_h)$ and the $S-$Order Runge-Kutta (RK) discretization method, defined as 
\begin{align*}
F_h(\bar{x}_h):=\bar{x}_h+h\sum_{k=1}^S b_kF_{es}(g_k),~g_k=\bar{x}_h+h\sum_{\ell=1}^{i-1}a_{ij}F_{es}(g_j),
\end{align*}
generate mappings $F_h$ that satisfy condition (a) in Definition \ref{regular_discretization}, \tcbb{where $a_{ij}\in\mathbb{R}$, $(i,j)\in S\times S$ defines the so called Runge-Kutta matrix \cite{numericalanalysisbook}}. However, since in the HAES \eqref{HAES}  the timer $\tau$ flows in the interval $[T_{\min},T_{\max}]$, any initial condition satisfying $\tau(0,0)<T_{\max}$ and located arbitrarily close to $T_{\max}$ could lead to discretized flows that generate updates of the form $\bar{\tau}_h^+=\bar{\tau}_h+hF_{\tau}(\bar{\tau}_h)>T_{\max}$, i.e., the solution leaves the flow set without hitting the jump set. To avoid this issue, we can consider a discretized jump set given by
\begin{equation}\label{jump_set2}
D_h:=D_{es}\cup\{\bar{x}_h:\bar{x}_h=F_h(y)\notin C,~y\in C_h\},
\end{equation}
which inflates the nominal jump set $D_{es}=\mathbb{R}^{n+m}\times\mathcal{T}_{D}\times\mathbb{T}^n$ in order to include the extra points outside of $C_{es}$ that may have been generated by discretized flows. The discretized flow set can then be defined as $C_h:=\mathbb{R}^{n+m}\times \mathcal{T}_C\times(\mathbb{T}^n+\rho(h))$, where $\rho\in\mathcal{K}_{\infty}$.

The following Proposition shows that the discretized HAES $\mathcal{H}_h:=\{C_h,F_h,D_h,\tcbb{G_{es}}\}$ retains the convergence properties of $\mathcal{H}_{a,\varepsilon}$ with respect to the set $\mathcal{A}$ given by \eqref{optimal}, provided the step size is sufficiently small.  The proof is a straightforward combination of \cite[Lemma 5.1]{SimulatorHybridSystems} with Proposition \ref{robustness_proposition} and the stability results of Theorems 1-4, and therefore it is omitted.
\begin{prop}\label{theorem_discretization}
Consider the discretized system $\mathcal{H}_h:=\{C_h,F_h,D_h,G_{es}\}$ with state $\bar{x}_h=(x_h,\tau_h,\mu_h)$, where $F_h$ is given either by the Forward Euler or S-Order RK discretization, and $D_h$ is given by \eqref{jump_set2}. Suppose that  Assumption \ref{assumption_freq} holds, as well as the Assumptions of Sections \ref{subsection:1}-\ref{subsection:4} for their respective HAES $\mathcal{H}_{a,\varepsilon}$. Then, the compact set $\mathcal{A}$ is SGPpAS as $(a,\varepsilon,h)\to0^+$ with $\mathcal{K}\mathcal{L}$ bound
\begin{equation}
|\bar{x}_h(\ell,j)|_{\mathcal{A}}\leq \beta_i(|\bar{x}_h(0,0)|_{\mathcal{A}},\ell h+j)+\nu,
\end{equation}
for all $(\ell,j)\in\text{dom}(\bar{x}_h)$, where $\ell\in\mathbb{Z}_{\geq0}$ is the index of the discretized flows, and $\beta_i$ is the $\mathcal{K}\mathcal{L}$ function generated by Theorem $i$, for $i\in\{1,2,3,4\}$.  
\end{prop}
Given that in Cases 2-4 the $\mathcal{K}\mathcal{L}$ bound $\beta_i$ is exponential, Proposition \ref{theorem_discretization} guarantees a geometric rate of convergence for their discretized HAES.  To guarantee completeness of solutions, the jump set can be slightly  modified as $G_h=(G_x^\top,T_{\min},\text{proj}_{\mathbb{T}^n}(\mu_h))^\top$, which simply sends $\mu_h$ back to $\mathbb{T}^n$ whenever there is a jump. 
\begin{figure}[t!]
  \centering
    \includegraphics[width=0.49\textwidth]{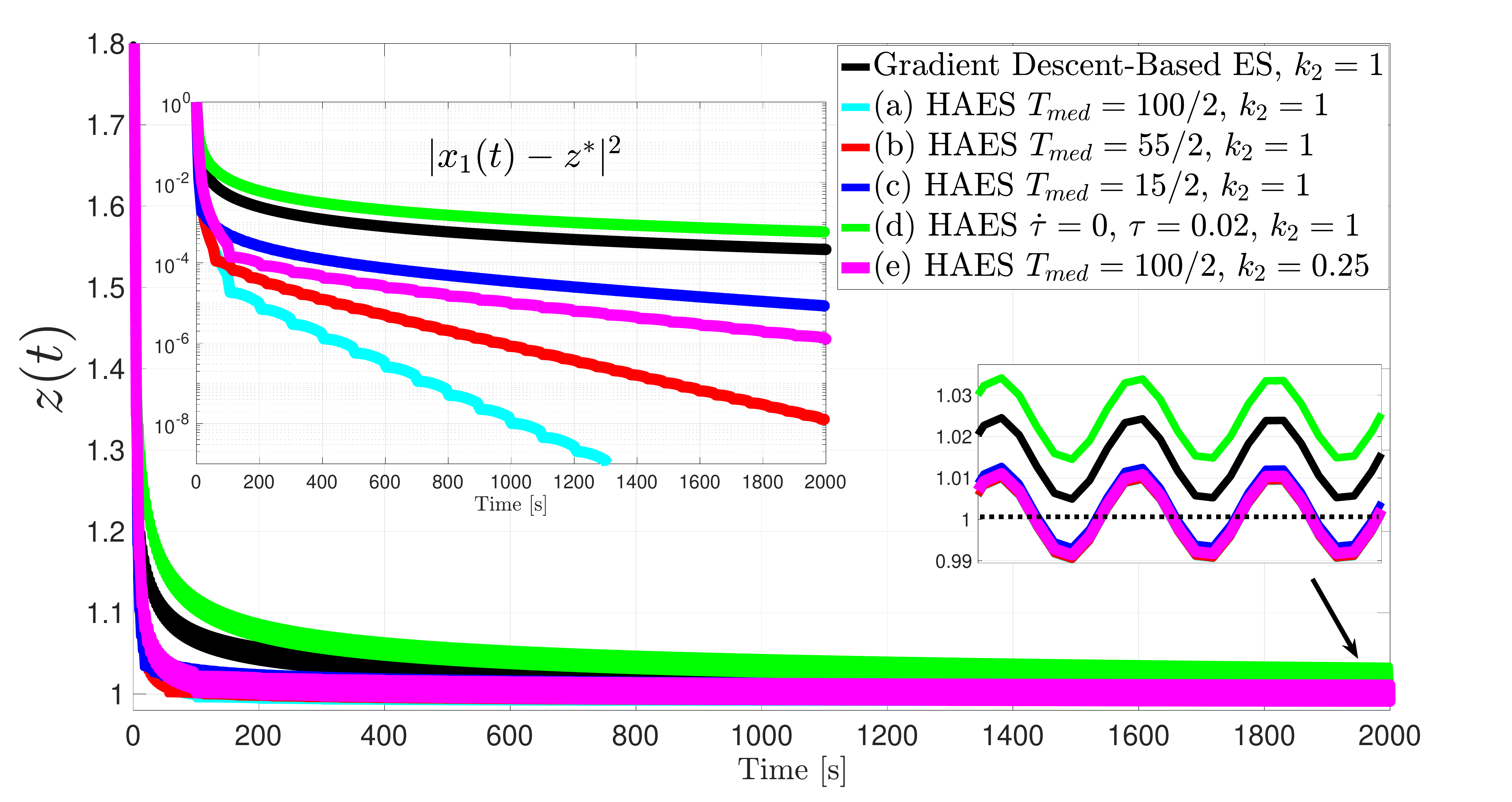}
      \caption{Solutions of the HAES and the classic gradient descent-based ES. The inset shows the different rates of convergence for the squared error $|x_1-z^*|^2$.}
      \label{figure1_convex}
\end{figure}
\section{Numerical Examples}
\label{section_numeric}
%
%
%
%
\subsection{Non-Strongly Convex Functions} 
\begin{figure}[t!]
  \centering
    \includegraphics[width=0.24\textwidth]{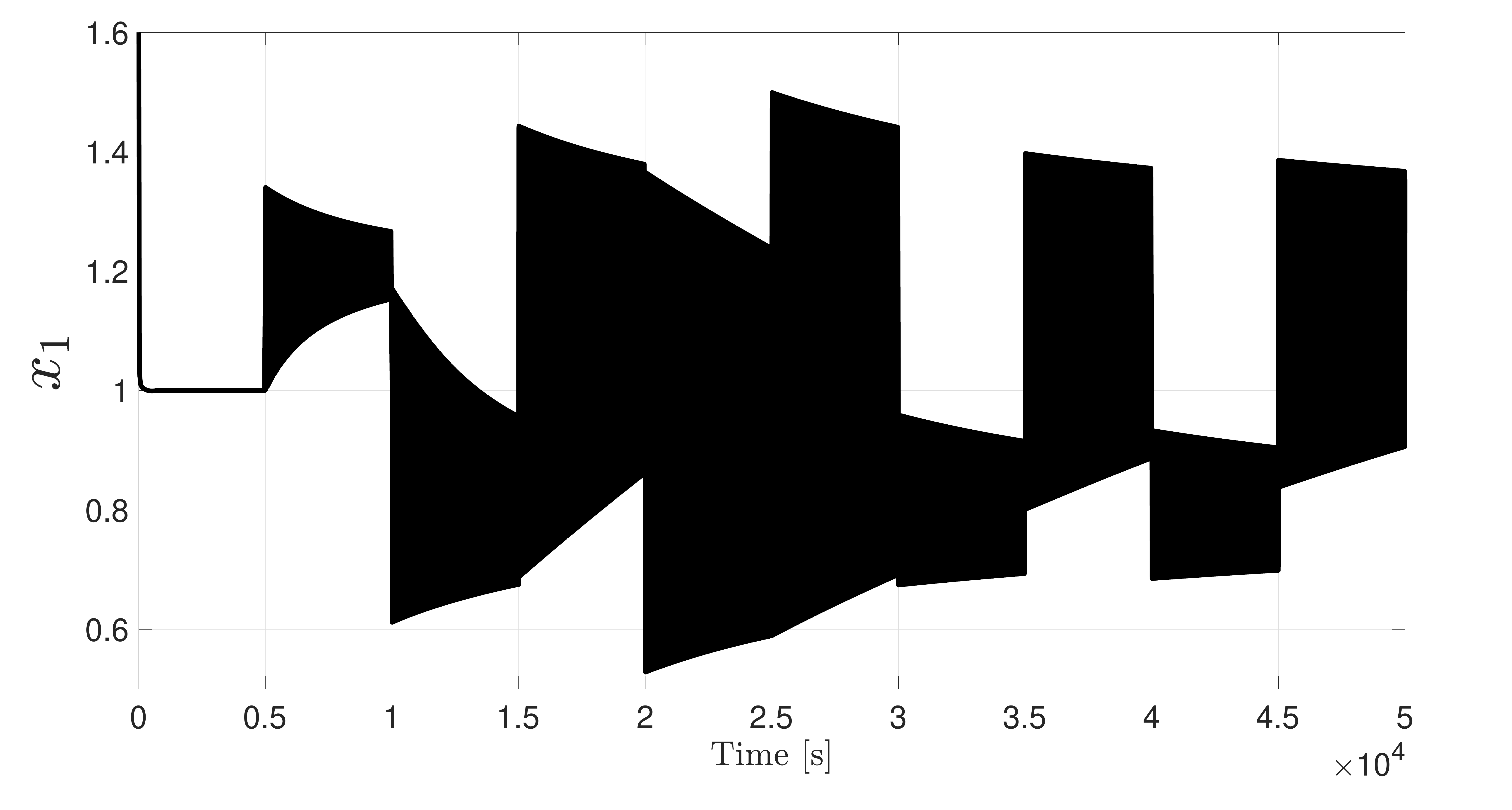}\includegraphics[width=0.24\textwidth]{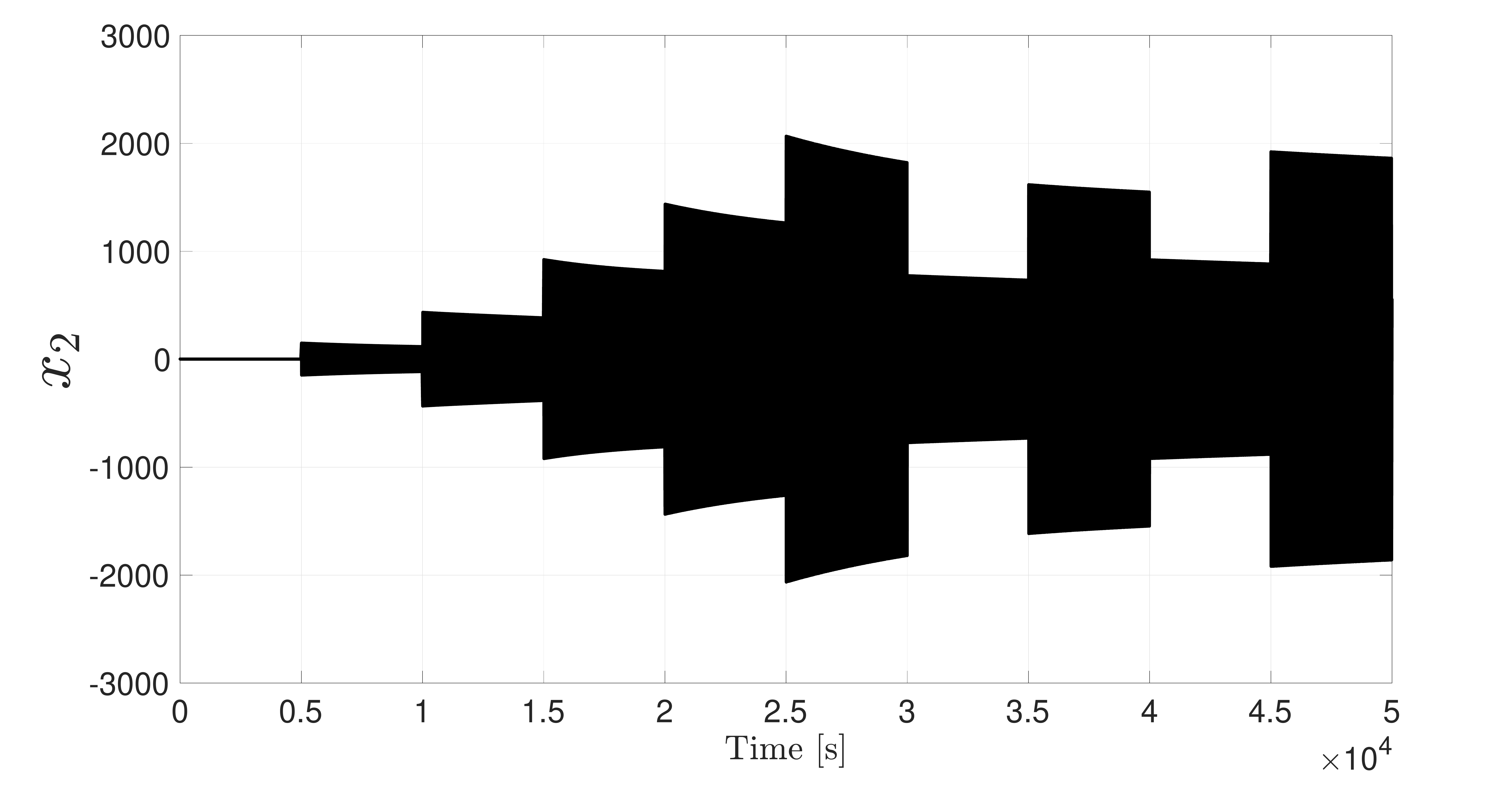}
    \includegraphics[width=0.24\textwidth]{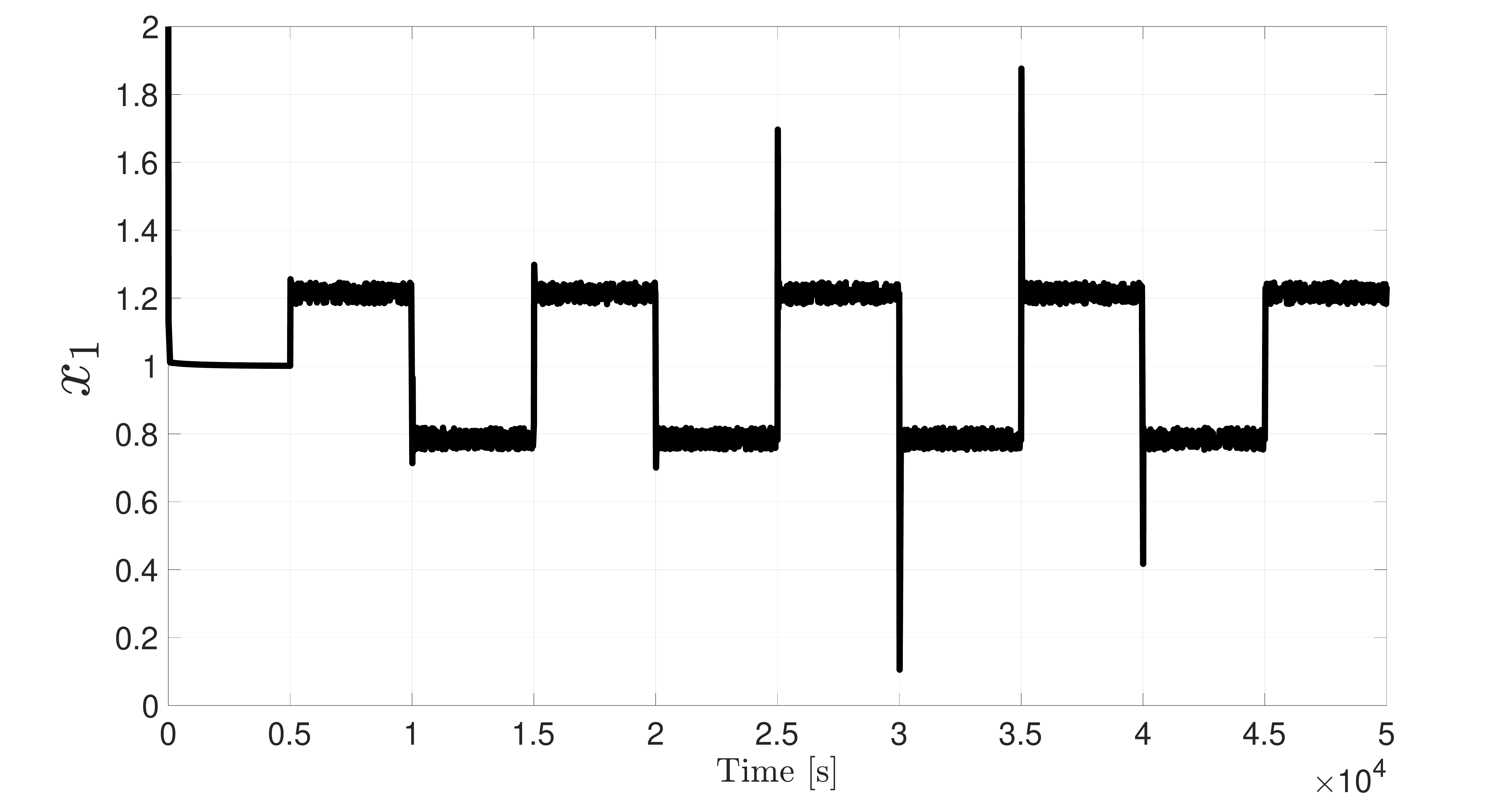}\includegraphics[width=0.24\textwidth]{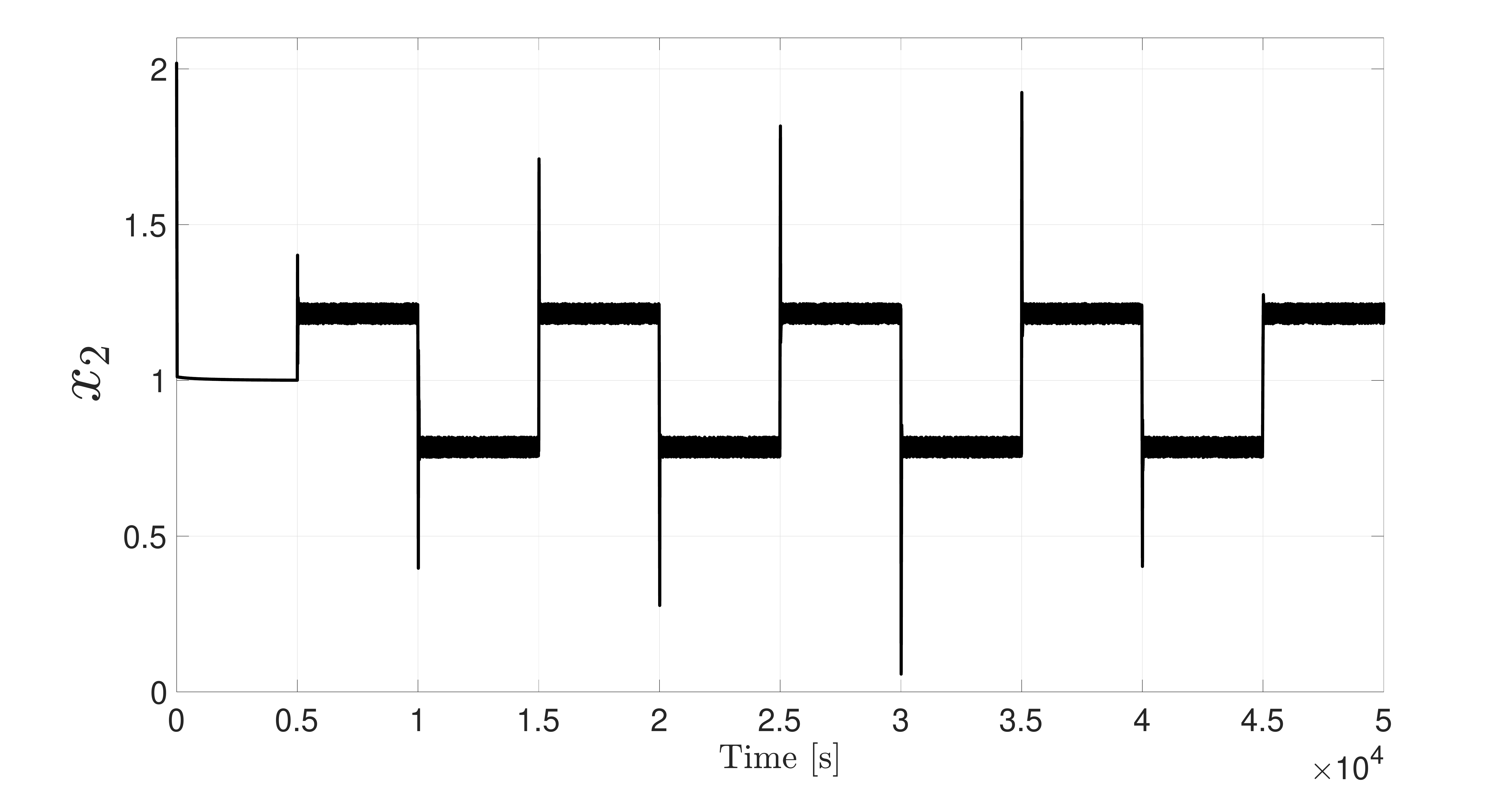}
      \caption{Solutions of the HAES under small disturbances with $T_{\text{med}}\to\infty$ (top) and with $T_{\text{med}}=T_{\text{max}}=25$ (bottom).}
      \label{figure2_convex}
\end{figure}

\begin{figure}[t!]
 \centering
   \includegraphics[width=0.49\textwidth]{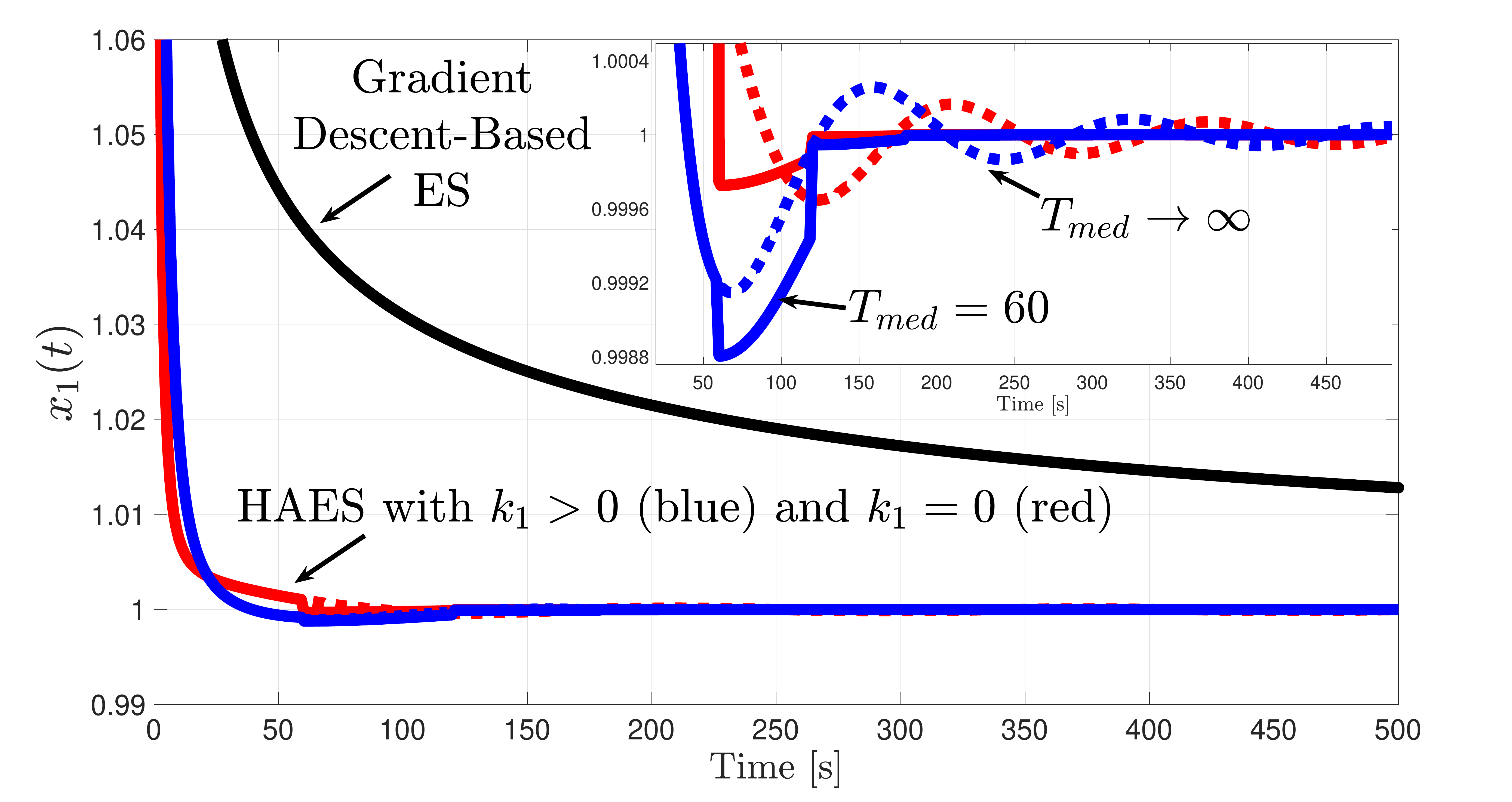}
      \caption{Effect of the Hessian damping and resseting in \eqref{maps1}. Red line corresponds to $k_1=0$ (no Hessian damping). Blue line corresponds to $k_1=10$. The inset shows how the frequent resettings remove the low frequency oscillations of the trajectories.}
      \label{NewFigure2}
\end{figure}

We first consider the function $\phi(z)=0.25(z-1)^{4}$, which is smooth, radially unbounded, convex, but not strongly convex. Since $\nabla^{2}\phi(z^*)=0$, this function does not satisfy the assumptions considered in \cite{Newton} and \cite{NewtonESC}. In order to achieve model-free optimization of $\phi$ with acceleration, we implement the HAES with data \eqref{maps1} and \eqref{sets_1} using a discretized HDS \eqref{discretized_HANDS} with discretized flow map obtained via 4-order Runge-Kutta method, discretized jump map \eqref{jump_set2}, and discretization step size $h=1\times 10^{-3}$. We consider the initial conditions $x_1(0,0)=2$, $x_2(0,0)=2$, $\tau(0,0)=0.01=T_{\min}$, and the parameters $k_1=0$, $k_2=1$, $a=0.01$, $\varepsilon=0.02$, $\kappa=2.54$. Figure \ref{figure1_convex} shows the evolution in time of five different solutions of the HAES, as well as a trajectory of the standard gradient descent-based ES dynamics using the same parameters $(a,\varepsilon)$. As shown in the plots, all solutions converge to a neighborhood of the optimal point $z^*$. However, as shown by the inset in the logarithmic scale, the rate of convergence is dramatically different for each solution. In particular, while solutions (a), (b), (c) and (e) converge to a small neighborhood of $z^*$ in approximately 40 seconds, the solution of the standard gradient descent-based ES algorithm requires almost 2000 seconds to reach the same neighborhood. This is consistent with the fact that solutions (a), (b), (c) and (e) exploit the acceleration property \eqref{inequality_algo1}. We also plotted solution (d), which keeps $\tau$ constant at $0.2$. In this case, the algorithm essentially approximates the time-invariant Heavy Ball ODE. As shown in Figure \ref{figure1_convex}, this dynamics generate the slowest rate of convergence, which illustrates the importance of the dynamic time-varying damping in \eqref{NesterovsODE}. The numerical results suggest that faster convergence is achieved by selecting large values of $T_{\text{med}}$, which is consistent with the bound \eqref{inequality_algo1}. However, as discussed in Section \ref{subsection:1}, in the limiting case when $T_{\text{med}}\to\infty$ the HAES behaves as the time-varying Nesterov's ODE with no restarting, which is highly sensitive to arbitrarily small disturbances. The top plots of Figure \ref{figure2_convex} illustrate this sensitivity, which emerges when $T_{\text{med}}>5\times 10^4$ and after adding a small disturbance $e(t)$ to the term $\phi(z)\tilde{\mu}$ in the flow map. The disturbance is a small squared periodic signal with frequency of $1\times 10^{-4}$ Hz and amplitude of $1\times 10^{-2}$. On the other hand, the robust stable behavior shown in the bottom plots correspond to the case when $T_{\text{med}}=T_{\max}=25$ and the same perturbation is added to the system. These simulations illustrate the importance of the restarting mechanism in ES algorithms with time-varying acceleration. 
Finally, we illustrate the effect of the Hessian damping term of \eqref{maps1}. Figure \ref{NewFigure2} shows the evolution of the state trajectory $x_1$ for $T_{med}=60$ and $T_{med}\to\infty$, and also for $k_1=0$ (no Hessian damping) and $k_1=10$ (with Hessian damping). In these simulations, we used $k_2=5$ (for both HAES and gradient descent-based ES), $\varepsilon=0.02$, $a=0.01$ and $h=1\times 10^{-4}$. As observed, the Hessian damping induced an initial slight improvement in the rate of convergence compared to the case $k_1=0$. In both cases ($k_1\geq0$) the transient performance of the HAES is significantly superior compared to the standard gradient-descent ES dynamics under the same gains and parameters ($k_2$, $a$, $\varepsilon$).
%
%
%
%
%
%
%
%
%

\subsection{Strongly Convex Functions}

\begin{figure}[t!]
  \centering
    \includegraphics[width=0.47\textwidth]{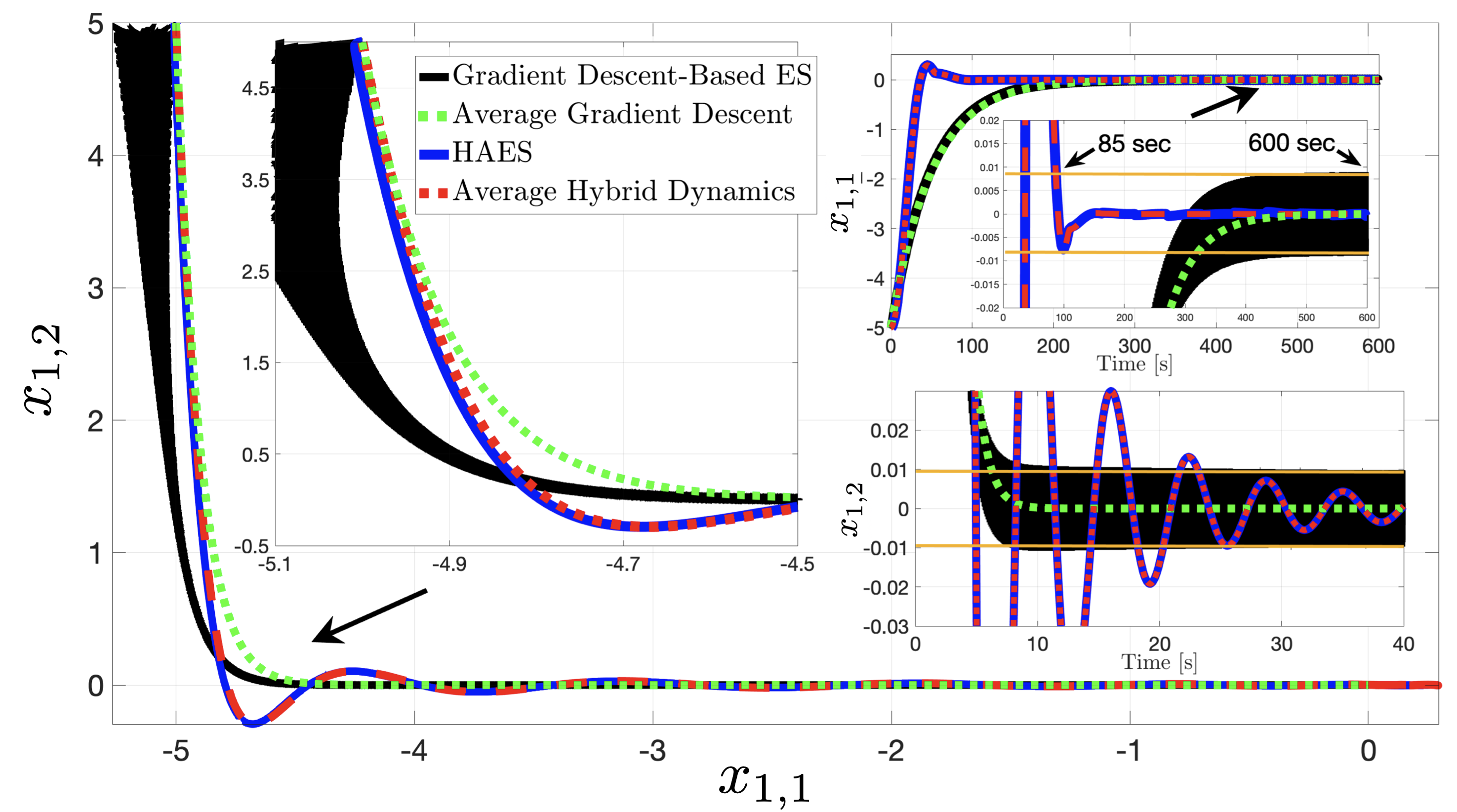}
      \caption{Comparison between trajectories generated by the HAES and the gradient descent-based ES. The inset shows the evolution in time of the optimizing states.}
      \label{Strong1}
\end{figure}

\tcbb{We now consider functions $\phi$ that satisfy Assumption \ref{assumption2}. First, we consider an ill-conditioned function in $\mathbb{R}^2$ given by $\phi(z)=\frac{1}{100}z_{1,1}^2+\frac{1}{2} z_{1,2}^2+10$. We simulate the standard gradient descent-based ES algorithm, as well as the HAES, both using the parameters $a=0.01$ and $\varepsilon=1\times10^{-3}$. The gain of the gradient descent-based ES algorithm is set to $k=1$, and the gain of the HAES is conservatively set to $k=0.25$. The restarting parameters $T_{\min}$ and $T_{\max}$ were set to $0.1$ and $27$, respectively. Figure \ref{Strong1} shows the resulting trajectories in the plane, as well as their evolution in time. It can be observed that both algorithms minimize $z_2$ at approximately the same speed, but the HAES minimizes $z_1$ approximately seven times faster. Moreover, the steady state oscillations of the HAES are substantially smaller. Figure \ref{Strong1} also shows the trajectories of the average systems related to both algorithms. The trajectory of the HAES is almost identical to the trajectory of its average hybrid dynamics. Finally, in order to illustrate the discussion of Section \ref{tuningstrong}, we also consider a multivariable ES problem of dimension $n=10$, where the cost function is given by $\phi(z)=\frac{1}{2}z^\top Qz+b^\top z+d$, with $d=10$, $b=[1,2,3,4,5,6,7,8,9,10]$, and $Q\in\mathbb{R}^{10\times 10}$ is a symmetric positive definite matrix generated randomly, which satisfies $L=\lambda_{\max}(Q)=33.01$, $\theta=\lambda_{\min}(Q)=1.757$. The parameters of the HAES are selected as $T_{\min}=0.01$, $a=0.01$, and $\varepsilon=1\times 10^{-3}$. Figure \ref{Level_sets} shows five trajectories of the HAES, as well as a trajectory of the gradient descent-based ES dynamics. With the exception of trajectory (a), which conservatively used $k=1/8L$, all other trajectories were generated by using the same gain $k=1/2L$. Trajectories (a) and (b) used the ``quasi-optimal'' restarting parameter $T_{\max}^*$ given by \eqref{quasi_optimal}. However, it can be observed that even for non-optimal values of $T_{\max}$, the hybrid dynamics significantly outperform the transient and steady state performance of the gradient descent-based ES algorithm, which decreases the sub-optimality measure at a rate of approximately $\theta/L\approx 0.0532$. On the other hand the trajectories generated by the HAES decrease the sub-optimality measure at a rate of approximately $\sqrt{\theta/L} \approx 0.23$. Both rates of decrease are indicated with dashed lines.}
Finally, we consider the cost function $\phi(z)=\frac{\|z\|^2}{4}$ with $z\in\mathbb{R}^2$, which has a Hessian matrix given by $\nabla^2 \phi(z)=0.5 I_2$. We compare the trajectories $x_1$ generated by the HAES \eqref{flows_generalHDS}-\eqref{jump_main} with mappings \eqref{maps2} and sets \eqref{sets_2} versus the trajectory generated by the standard gradient descent-based ES dynamics \cite{tan06Auto}. Figure \ref{Level_sets} shows both trajectories evolving over the level sets of the cost function. It can be observed that the HAES exhibits significant less oscillations compared to the standard gradient descent-based ES algorithm using the same parameters $(k_2,a,\varepsilon)$.
\begin{figure}[t!]
  \centering
    \includegraphics[width=0.49\textwidth]{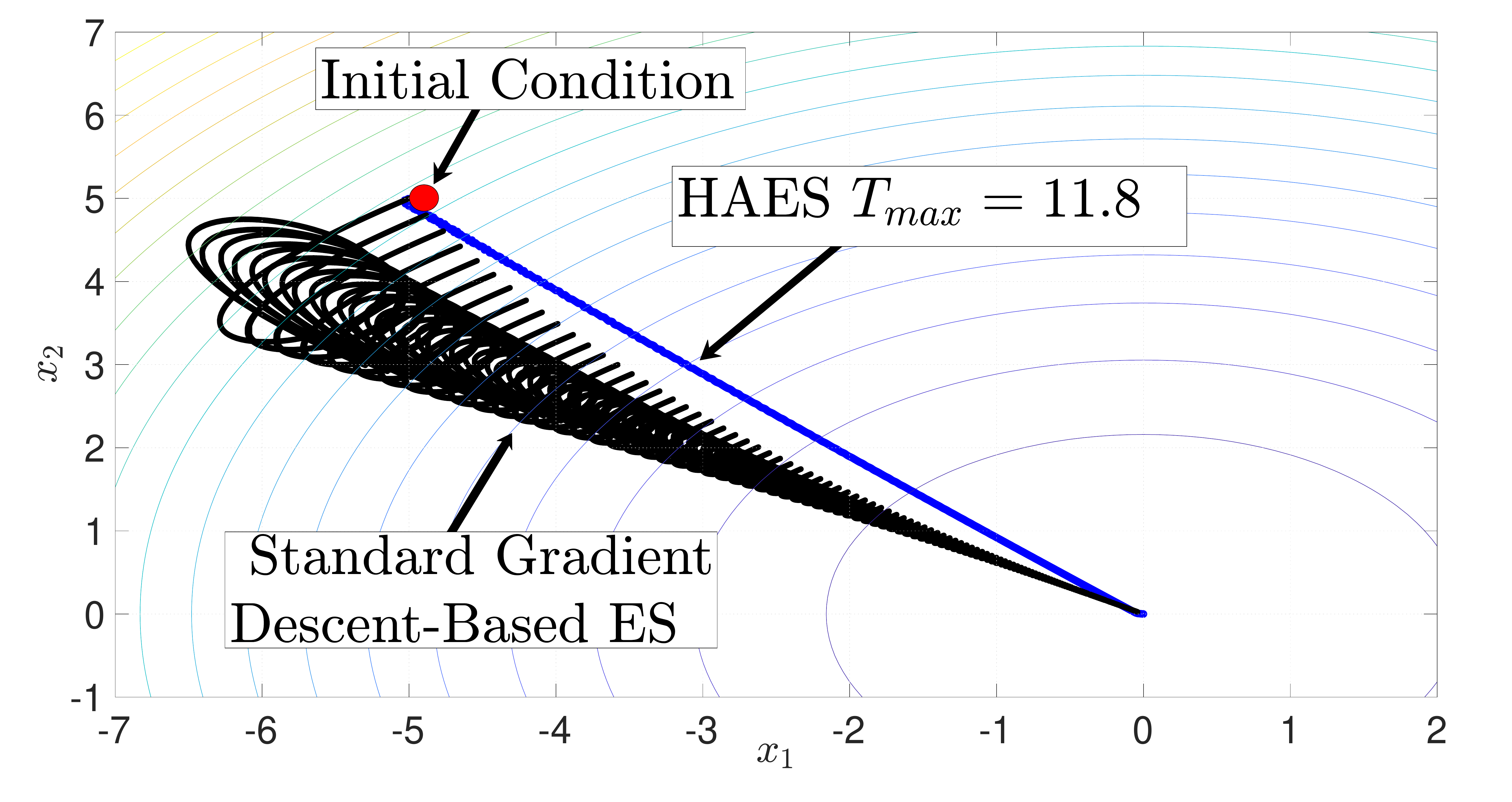}
      \caption{Evolution in $\mathbb{R}^2$ of the solutions of the derivative-free hybrid dynamics, and the solutions of the classic gradient descent-based extremum seeking algorithm.}
      \label{Level_sets}
\end{figure}
\section{Analysis: Part 1 - Averaging Theory}
\label{averaging_section}
In order to prove the main results of this paper, we first develop some auxiliary stability results for a class of singularly perturbed HDS \cite{Wang:12_Automatica,averagingTeel} that fits the structure of our algorithms. In particular, we consider HDS with states $(\varphi,\chi)\in\mathbb{R}^{n_1}\times\mathbb{R}^{n_2}$, and continuous-time dynamics parameterized by two different constants $\varepsilon>0$ and $\delta>0$. The parameter $\varepsilon>0$ induces a multi-time scale behavior in  the flow map. The constant $\delta$ parametrizes the stability properties of the slow dynamics. The singularly perturbed HDS is modeled as: 
\begin{subequations}\label{SPHDS}
\begin{align}
&\dot{\varphi}=f^{\delta}_{\varphi}(\varphi,\chi),~\dot{\chi}=\dfrac{1}{\varepsilon}f_{\chi}(\varphi,\chi),~~(\varphi,\chi)\in C\times \Psi\\
&\varphi^+\in G_{\varphi}(\varphi,\chi),~~~\chi^+= \chi,~~~~(\varphi,\chi)\in D\times \Psi,
\end{align}
\end{subequations}
where $\varphi\in\mathbb{R}^{n_1}$, $\chi\in\mathbb{R}^{n_2}$, $C,D\subset\mathbb{R}^{n_1}$, $\Psi\subset\mathbb{R}^{n_2}$, $f^{\delta}_{\varphi}:\mathbb{R}^{n_1}\times\mathbb{R}^{n_2}\to\mathbb{R}^{n_1}$, $f_{\chi}:\mathbb{R}^{n_1}\times\mathbb{R}^{n_2}\to\mathbb{R}^{n_2}$, $G_{\varphi}:\mathbb{R}^{n_1}\times\mathbb{R}^{n_2}\rightrightarrows\mathbb{R}^{n_1}$ is a set-valued mapping, and $\delta,\varepsilon\in\mathbb{R}_{>0}$. For the sake of generality we will allow set-valued jump maps $G_{\varphi}$, as well as stability results that are local with respect to some basin of attraction.  We make the following regularity assumption on system \eqref{SPHDS}.
\begin{figure}[t!]
  \centering
    \includegraphics[width=0.5\textwidth]{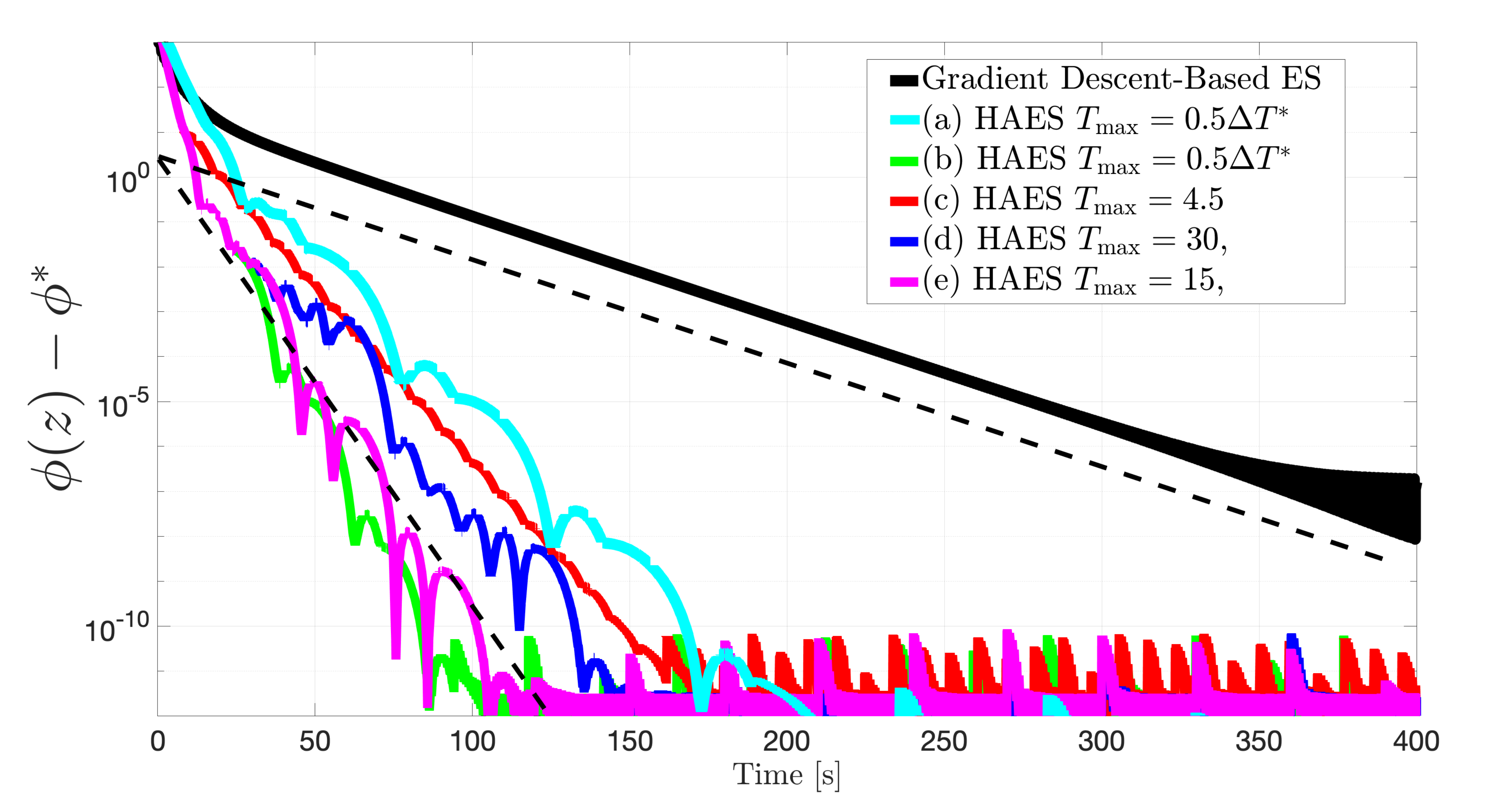}
      \caption{\tcbb{Decrease of the sub-optimality measure under gradient descent-based ES and the HAES, for a multivariable problem with $n=10$.}}
      \label{Level_sets}
\end{figure}
\begin{assumption}\label{assump1}
For each $\delta>0$ the following holds: The sets $C$ and $D$ are closed, the functions $f^{\delta}_{\varphi}$ and $f_{\chi}$ are continuous for all $(\varphi,\chi)\in C\times\Psi$, the set-valued mapping $G_{\varphi}$ is outer-semicontinuous \cite[Def. 5.9]{Goebel:12} and locally bounded \cite[Def. 5.14]{Goebel:12}, $\Psi$ is compact, and for each $(\varphi,\chi)\in D\times\Psi$ the set $G_{\varphi}(\varphi,\chi)$ is not empty.  
\end{assumption}

To analyze the HDS \eqref{SPHDS}, the hybrid dynamics are expressed in the time variables $(\tau,j)$, where $\tau=t/\varepsilon$, i.e., 
\begin{subequations}\label{SPHDS2}
\begin{align}
&
\dfrac{d\varphi}{d\tau}=\varepsilon f^{\delta}_{\varphi}(\varphi,\chi),~\dfrac{d\chi}{d\tau}=f_{\chi}(\varphi,\chi),~~(\varphi,\chi)\in C\times \Psi\\
&
\varphi^+\in G_{\varphi}(\varphi,\chi),~~~\chi^+=\chi,~~~~~~~~~~~(\varphi,\chi)\in D\times \Psi.
\end{align}
\end{subequations}
For this system we define the \emph{boundary layer dynamics}.
\begin{definition}\label{defBDL}
The boundary layer dynamics of the hybrid system \eqref{SPHDS2} are given by $(\varphi_{bl},{\chi}_{bl})\in C\times \Psi$, $\frac{d\varphi_{bl}}{d\tau}=0,~\frac{d{\chi}_{bl}}{d\tau}=f_{\chi}(\varphi_{bl},{\chi}_{bl})$, which ignores the jumps, and ``freezes'' $\varphi$ by setting $\varepsilon=0$.
\end{definition}
Similar to existing results in singular perturbation and averaging theory, e.g., \cite{Wang:12_Automatica}, our goal is to establish stability properties for the singularly perturbed hybrid system \eqref{SPHDS} based on a simplified \emph{average} system obtained by averaging the dynamics of $\varphi$ along the solutions of ${\chi}$. To do this we need the following assumption.
\begin{assumption}\label{assump2}
For each $\delta>0$ there exists a continuous function $f^{\delta}_A:\mathbb{R}^{n_1}\to\mathbb{R}$ such that for each compact set $K\subset C\times \Psi$ there exists a class-$\mathcal{L}$ function $\sigma_{K,\delta}$ such that, for each $L>0$, each $\varphi\in C\cap K$, and each function ${\chi}_{bl}:[0,L]\to\Psi$ satisfying $\dot{{\chi}}_{bl}=f_{\chi}(\varphi,{\chi}_{bl})$, the following holds:
\begin{equation}\label{bound_average}
\left|\frac{1}{L}\int^L_{0}\Big(f_{\varphi}^{\delta}(\varphi,{\chi}_{bl}(s))-f^{\delta}_{A}(\varphi)\Big)ds\right|\leq \sigma_{K,\delta}(L).
\end{equation}
\end{assumption}
Using the mapping $f^{\delta}_A$ defined in \eqref{bound_average} we now define the average hybrid system of \eqref{SPHDS}:
\begin{definition}\label{definition_average_system}
The average hybrid system $\mathcal{H}_{\delta}^A$ of the HDS \eqref{SPHDS} with boundary layer dynamics specified in Definition \ref{defBDL} has a state $y\in\mathbb{R}^{n_1}$, and is given by
\begin{align}\label{SPHDS_average}
\dot{y}=f^{\delta}_A(y),~~y\in C,~~~~~~~~y^+\in G_A(y),~~y\in D,
\end{align}
where $G_A:=\{v_1\in\mathbb{R}^{n_1}:(v_1,v_2)\in G_{\varphi}(\varphi,\chi),(\chi,v_2)\in\Psi\times\mathbb{R}^{n_2}\}$.
\end{definition}
Finally, we assume that the average system \eqref{SPHDS_average} satisfies the following semi-global practical pre-stability property with respect to a compact set $\mathcal{A}_{\varphi}$.
\begin{assumption}\label{assumption_stability}
There exists a nonmepty compact set $\mathcal{A}_{\varphi}\subset\mathbb{R}^{n_1}$, an open set $\mathcal{B}_{\mathcal{A}_{\varphi}}\supset\mathcal{A}_{\varphi}$, and a class $\mathcal{K}\mathcal{L}$ function $\beta$ such that for each proper indicator\footnote{A proper indicator of $\mathcal{A}_{\varphi}$ on $\mathcal{B}_{\mathcal{A}_{\varphi}}$ is a continuous function $\omega:\mathcal{B}_{\mathcal{A}_{\varphi}}\to\mathbb{R}_{\geq0}$ satisfying $\omega(\varphi)=0$ if and only if $\varphi\in\mathcal{A}_{\varphi}$, and such that $\omega(\varphi_i)\to\infty$ when $i\to\infty$ if either $|\varphi_i|\to\infty$ or the sequence $\{\varphi_i\}_{i=1}^n$ approaches the boundary of $\mathcal{B}_{\mathcal{A}_{\varphi}}$.} $\omega(\cdot)$ for $\mathcal{A}_{\varphi}$ on $\mathcal{B}_{\mathcal{A}_{\varphi}}$, each  compact set $K_0\subset\mathcal{B}_{\mathcal{A}_{\varphi}}$, and each $\nu>0$, there exists a $\delta^*>0$ such that for all $\delta\in(0,\delta^*)$, all solutions of \eqref{SPHDS_average} with $y(0,0)\in K_0$ satisfy the bound:
\begin{equation}\label{KL_average_system}
\omega(y(t,j))\leq\beta (\omega(y(0,0)),t+j)+\nu. 
\end{equation}
 for all $(t,j)\in\text{dom}(y)$. 
\end{assumption}
Using Assumptions \ref{assump1}, \ref{assump2}, and \ref{assumption_stability} the following two results are obtained. Proposition \ref{closeness_on_compact_time_domain} is a straightforward extension of \cite[Thm. 1]{Wang:12_Automatica} that concerns closeness on compact time domains of the $\varphi$-component of the solutions of \eqref{SPHDS} to the solutions $y$ of \eqref{SPHDS_average}. On the other hand, Theorem \ref{main_theorem_averaging} links the stability properties of system \eqref{SPHDS} to the stability properties of the average system \eqref{SPHDS_average}.
\begin{prop}\label{closeness_on_compact_time_domain}
Suppose that the HDS \eqref{SPHDS} satisfies Assumptions \ref{assump1}, \ref{assump2}, and \ref{assumption_stability}. Let $K_0\subset\mathcal{B}_\mathcal{A}$ and let $\delta^*>0$ be such that  $\forall~\delta\in(0,\delta^*)$ all solutions of \eqref{SPHDS_average} with $y(0,0)\in K_0$ do not have finite escape times. Then, for each $\delta\in(0,\delta^*)$, each $\rho>0$, and any strictly positive real numbers $T,J$ there exists $\varepsilon^*>0$ such that for each $\varepsilon\in(0,\varepsilon^*]$ and each solution $\varphi$ to system \eqref{SPHDS} with $\varphi(0,0)\in K_0$, there exists some solution $y$ to the average system \eqref{SPHDS_average}  with $y(0,0)\in K_0$ such that $\varphi$ and $y$ are $(T,J,\rho)$-close in the sense of \cite[Def. 5.23]{Goebel:12}.
\end{prop}
\noindent
\textbf{Proof}: Since $\forall~\delta>0$ the singularly perturbed HDS \eqref{SPHDS} satisfies all the assumptions needed to apply \cite[Thm. 1]{Wang:12_Automatica}, it only needs to be shown that there exists a $\delta^*>0$ such that the system has no finite escape times from $K_0$. Indeed, by Assumption \ref{assumption_stability} for each compact set of initial conditions $K_0$ there exists $\delta^*>0$ such that for all $\delta\in(0,\delta^*)$ all the solutions satisfy \eqref{KL_average_system}, which precludes finite escape times. \null \hfill \null $\blacksquare$

%
\begin{thm}\label{main_theorem_averaging}
Suppose that the HDS \eqref{SPHDS} satisfies Assumptions \ref{assump1}, \ref{assump2}, and \ref{assumption_stability}. Then, for each proper indicator $\omega$ for $\mathcal{A}_{\varphi}$ on $\mathcal{B}_{\mathcal{A}_{\varphi}}$, each compact set $K_0\subset\mathcal{B}_{\mathcal{A}_{\varphi}}$ and each $\nu>0$ there exists $\delta^*>0$ such that for each $\delta\in(0,\delta^*)$ there exists $\varepsilon^*>0$ such that for all $\varepsilon\in(0,\varepsilon^*)$ all solutions of \eqref{SPHDS} with $\varphi(0,0)\in K_0$ satisfy $\omega(\varphi(t,j))\leq\beta (\omega(\varphi(0,0)),t+j)+\nu$ for all $(t,j)\in\text{dom}(\varphi)$. 
\end{thm}
\textbf{Proof}: The proof is similar to the proofs of \cite[Thm. 2]{Wang:12_Automatica} and \cite[Thm. 2]{averagingTeel}. Let $K_0\subset\mathcal{B}_{\mathcal{A}_{\varphi}}$ and $\nu>0$ be given. Let $\omega:\mathcal{B}_{\mathcal{A}_{\varphi}}\to\mathbb{R}_{\geq0}$ be a proper indicator for $\mathcal{A}_{\varphi}$ with respect to $\mathcal{B}_{\mathcal{A}_{\varphi}}$. Define the set $K_1:=\left\{\varphi\in \mathcal{B}_{\mathcal{A}_{\varphi}}:~\omega(\varphi)\leq\beta\left(\max_{y\in K_0}\omega(y),0\right)+1\right\}$, and

\begin{equation}\label{big_K}
K:=K_1\cup G_A(K_1\cap D).
\end{equation}

Since $K_1$ is compact, and $G_A$ is outer semicontinuous and locally bounded, the set $K$ is compact. Moreover $K\subset\mathcal{B}_{\mathcal{A}_{\varphi}}$ since $\omega$ is a proper indicator and $G_A$ is an OSC mapping that maps $\mathcal{B}_{\mathcal{A}_{\varphi}}\cap D$ to $\mathcal{B}_{\mathcal{A}_{\varphi}}$.  Let $\epsilon_1>0$ be such that, for all $\varphi\in K$, all $y\in K+\epsilon_1\mathbb{B}$ with $|\varphi-y|\leq \epsilon_{1}$, and all $s\geq0$, the following holds:

\begin{align}\label{inequa1}
\omega(\varphi)\leq \omega(y)+\frac{\nu}{3},~~~~~\beta(\omega(y),s)\leq \beta(\omega(\varphi),s)+\frac{\nu}{3}
\end{align}

Such $\epsilon_1^*>0$ always exists given that $\beta\in\mathcal{K}\mathcal{L}$ and that $\beta,\omega$ are continuous functions. Using Proposition \ref{robustness_proposition} in the Appendix \ref{AppendixA}, there exists a $\delta^*>0$ such that for all $\delta\in(0,\delta^*)$ there exists a $\rho^*\in(0,\epsilon_1)$ such that for all $\rho\in(0,\rho^*)$ all solutions $y$ of the $\rho$-inflation of system \eqref{SPHDS_average} with $y(0,0)\in K$ satisfy for all $(t,j)\in\text{dom}(y)$ the following bound:
\begin{equation}\label{bound_KL_ave_infla}
\omega(y(t,j))\leq\beta (\omega(y(0,0)),t+j)+\frac{\nu}{3}. 
\end{equation}

Let $\mu\geq0$ and consider the extended hybrid dynamical system, constructed from \eqref{SPHDS}, with auxiliary state $\eta\in\mathbb{R}^{n_1}$ and $K$-restricted flow and jump set, given by:
\begin{subequations}\label{restricted}
\begin{align}
&\left.\begin{array}{l}
\dot{\varphi}=f^{\delta}_{\varphi}(\varphi,\chi),~~~\dot{\chi}=\frac{1}{\varepsilon}f_{\chi}(\varphi,\chi)\\
\dot{\eta}=\frac{1}{\varepsilon}\left[f_{\varphi}^{\delta}(\varphi,\chi)-f^{\delta}_A(\varphi)-\mu\eta\right]
\end{array}\right\},\begin{array}{l}(\varphi,\chi,\eta)\\
\in (C\cap K)\times \Psi\times\mathbb{R}^n\end{array}\\
&\left.\begin{array}{l}
\varphi^+\in G_{\varphi}(\varphi,\chi)\\
{\chi}^+= {\chi},~\eta^+=0
\end{array}\right\},~(\varphi,{\chi},\eta)\in (D\cap K)\times \Psi\times\mathbb{R}^n.
\end{align}
\end{subequations}
Since for each $\delta>0$ all the assumptions needed to apply \cite[Lemma 4]{Wang:12_Automatica} are satisfied, the next Lemma follows directly by \cite[Lemma 4]{Wang:12_Automatica}.
\begin{lem}\label{lemma_aux1}
Suppose that the HDS \eqref{SPHDS} satisfies Assumptions \ref{assump1} and \ref{assump2}. Then, for each $\delta\in(0,\delta_{K}^*)$ and each $\rho>0$ there exists $\varepsilon^*,\lambda\in\mathbb{R}_{>0}$ such that, for all $\varepsilon\in(0,\varepsilon^*]$, each solution $(\varphi,\chi,\eta)$ of system \eqref{restricted} with $\eta(0,0)=0$ satisfies $\lambda|\eta(t,j)|\leq \rho$, for all $(t,j)\in\text{dom}(\varphi,\chi)$. \QEDB
\end{lem}
%
%

\noindent
Let $K$, $\delta$, and $\rho$, generate $(\epsilon_2,\lambda)$ via Lemma \ref{lemma_aux1}. Let $\varepsilon^*:=\min\{\rho,\epsilon_2,\lambda\}$ and let $\varepsilon\in(0,\varepsilon^*)$. For each solution $(\varphi,\chi,\eta)$ of \eqref{restricted} let us define $y:=\varphi-\varepsilon\eta$. Since $\eta^+=0$, we obtain $\dot{y}=\dot{\varphi}-\varepsilon\dot{\eta}=f^\delta_A(\varphi)+\lambda\eta$ and $y^+=\varphi^+\in G_A(\varphi)$ with
\begin{align*}
G_A(\varphi):= \{v_1\in\mathbb{R}^{n_1}:(v_1,v_2)\in G_{\varphi}(\varphi,\chi),(\chi,v_2)\in\Psi\times\mathbb{R}^{n_2}\}.
\end{align*}
Since $\varphi=y+\varepsilon\eta$, we have that
\begin{align*}
\dot{y}&=f^\delta_A(y+\varepsilon\eta)+\lambda\eta,~~~~~y+\varepsilon\eta\in C,\\
y^+&\in G_A(y+\varepsilon\eta),~~~~~~~~~~~~y+\varepsilon\eta\in D,
\end{align*}
and by the choice of $\varepsilon^*$ above we have that $\dot{y}\in f^\delta_A(y+\rho\mathbb{B})+\rho\mathbb{B}$, when $y\in C_{\rho}$, and $y^+\in G_A(y+\rho\mathbb{B})$, when $y\in D_{\rho}$, where the sets $C_\rho$ and $D_\rho$ correspond to the $\rho$-inflations constructed as in \eqref{C_inflated} and \eqref{D_inflated}. Therefore, for each solution $(\varphi,\chi,\eta)$ of system \eqref{restricted} with $(\varphi(0,0),\chi(0,0))\in K_0\times\Psi$ and $\eta(0,0)=0$, the trajectory $y:=\varphi-\varepsilon\eta$ is a solution to the $\rho$-inflation \eqref{perturbedHDS} of the average system \eqref{SPHDS_average}, and since $\delta\in(0,\delta^*)$ the trajectory $y$ also satisfies the bound \eqref{bound_KL_ave_infla}. Since $\rho\leq  \epsilon^*_1$, the inequalities \eqref{inequa1} hold, and all solutions $(\varphi,\chi)$ to the system \eqref{restricted} with $(\varphi(0,0),\chi(0,0))\in K_0\times \Psi$ satisfy for all $(t,j)\in\text{dom}(\varphi,\chi)$ the following bounds:
\begin{align}
\omega(\varphi(t,j))&\leq \omega(y(t,j))+\frac{\nu}{3}\leq\beta\left(\omega(y(0,0)),t+j\right)+\frac{2\nu}{3}\notag\\
&\leq\beta\left(\omega(\varphi(0,0)),t+j\right)+\nu\label{klbound}.
\end{align}
Since $\nu\in(0,1)$, by the inequality \eqref{klbound}, each solution of \eqref{restricted} with $\varphi(0,0)\in K_0$ remains in the compact set  $K_v:=\left\{\varphi\in\mathbb{R}^{n_1}:~\omega(\varphi)\leq \beta(\max_{\bar{\varphi}\in K_0}\omega(\bar{\varphi}),0)+v\right\}$,  which is contained in the interior of the set $K$.

We now use the properties of the solutions of the $K$-restricted system \eqref{restricted} to derive conclusions about the solutions of the original HDS \eqref{SPHDS}. Indeed, since $K_0\subset K$, a solution of \eqref{SPHDS} with $(\varphi(0,0),\chi(0,0))\in K_0\times \Psi$ must agree with a solution of \eqref{restricted} for all $(t,j)\in\text{dom}(\varphi,\chi)$ such that $\varphi(t,j)\in K$. However, using the definition of $K$ in \eqref{big_K} and the $\mathcal{K}\mathcal{L}$ bound \eqref{klbound} we have that all solutions of \eqref{SPHDS} with $(\varphi(0,0),\chi(0,0))\in K_0\times \Psi$ remain in the set $\left(K_v\cup G_A(K_v\cap D)\right)\times \Psi\subset K\times\Psi$. This implies that inequality \eqref{klbound} holds for all $(t,j)\in\text{dom}(\varphi,\chi)$, which establishes the result. \null \hfill \null $\blacksquare$
\section{Analysis: Part 2 - Algorithmic Stability}
\label{sec_analysis}
In this section we use Theorem \ref{main_theorem_averaging} to prove Theorems \ref{thm1a}-\ref{thm4}.
%
%
%
%
%
%
In particular, we show that all HAES can be written as a singularly perturbed system of the form \eqref{SPHDS} with $\delta=a$, and that all the assumptions needed to apply Theorem \ref{main_theorem_averaging} hold.  Indeed, by construction, it can be seen that for small values of $\varepsilon>0$ 
%
the HDS \eqref{flows_main_equations}-\eqref{main_jump_dynamics} is a singularly-perturbed HDS of the form  \eqref{SPHDS} with $\varphi=(x,\tau)$, $\chi=\mu$, $C=\mathbb{R}^{n+m}\times\mathcal{T}_C$, $D=\mathbb{R}^{n+m}\times\mathcal{T}_D$, and $\Psi=\mathbb{T}^n$. By construction of the dynamics, Assumption \ref{assump1} is satisfied since for each $a>0$ all the mappings $F_x$ and $G_x$ are continuous in $C\times \Psi$ and $D\times\Psi$, respectively, and the sets $\mathcal{T}_C$ and $\mathcal{T}_D$ are closed. 
\subsection{Average Hybrid Systems}
\label{section_average_systems}
We now show that all the HAES satisfy Assumption \ref{assump2}. The following Lemma, which relies on Assumption \ref{assumption_freq}, will be instrumental for our results. The proof is presented in the Appendix C.
\begin{lemma}\label{lemma_dither}
Suppose that Assumption \ref{assumption_freq} holds. Then, there exists a $\tilde{T}>0$ such that every solution $\mu:\mathbb{R}_{\geq0}\to\mathbb{R}^{2n}$ of $\dot{\mu}=R\mu$ with $\mu(0)\in\mathbb{T}^n$ satisfies $\frac{1}{N\tilde{T}}\int_{0}^{N\tilde{T}}\tilde{\mu}(\tau)\tilde{\mu}(\tau)^\top d\tau=\frac{1}{2}I_{n}$, and $\int_{0}^{N\tilde{T}}\tilde{\mu}(\tau)d\tau=\mathbf{0}_n$, for any $N\in\mathbb{Z}_{>0}$, where $\tilde{\mu}$ is defined in \eqref{input_1}.
\end{lemma}
Since the cost function $\phi$ is at least twice continuously differentiable, the 2nd-Taylor expansion of  $\phi(x_1+a\tilde{\mu})$ around $x_1$ is well defined and given by $\phi(x_1+a\tilde{\mu})=\phi(x_1)+a\tilde{\mu}^\top \nabla \phi(x_1)+\mathcal{O}(a^2)$. Substituting in \eqref{maps1}, \eqref{maps2}, \eqref{map3}, and \eqref{map4}, and using the fact that $|\tilde{\mu}|\leq 1$ and introducing the function $\Upsilon(x_1,\tilde{\mu},a):=\frac{\phi(x_1)}{a}\tilde{\mu}+\tilde{\mu}\tilde{\mu}^\top \nabla \phi(x_1)$, for each of the Cases 1-4, we obtain mappings $F_x:=(F_{x_1}^\top,F_{x_2}^\top)^\top$ with components: 
\begin{enumerate}[(a)]
\item For Cases 1 and 2: $F_{x_1}=\dfrac{2}{\tau}(x_2-x_1)-2k_1\Upsilon(x_1,\tilde{\mu},a)$ and 
$F_{x_2}=-4k_2\tau\Upsilon(x_1,\tilde{\mu},a)+\mathcal{O}(a)$.
\item For Case 3: $F_{x_1}=-2\Upsilon(x_1,\tilde{\mu},a)-k A^\top x_2+\mathcal{O}(a)$
\item For Case 4: $F_{x_1}=-2\Upsilon(x_1,\tilde{\mu},a)-k\sum_{j=1}^mH_j(x)A_j+\mathcal{O}(a).$
%
%
\end{enumerate}
For each Case 1-4, consider the mapping
\begin{equation}\label{integral1}
f^a_A(x,\tau)=[f^{a\top}_{A,1}, f^{a\top}_{A,2},f^a_{A,3}]^\top,
\end{equation}
defined as follows. For Cases 1 and 2: $f^a_{A,1}=\dfrac{2}{\tau}(x_2-x_1)-k_1\nabla\phi(x_1)+\mathcal{O}(a)$, $f^a_{A,2}=-2k_2\tau \nabla \phi(x_1)+\mathcal{O}(a)$, $f^a_{A,3}=\{0.5,1\}$. For Case 3: $f^a_{A,1}=-\nabla \phi(x_1)-kA^\top x_2+\mathcal{O}(a)$, $f^a_{A,2}=Ax_1-b$, and $f^a_{A,3}=0$.  For Case 4: $f^a_{A,1}=-\nabla \phi(x_1)-k\sum_{j=1}^mH_j(x)A_j+\mathcal{O}(a)$, $f^a_{A,2}=\sum_{j=1}^n\left(H_j(x)-x_{2,j}\right)e_j$, $f^a_{A,3}=0$. Under these definitions, the following lemma shows that the mapping \eqref{integral1} satisfies the conditions of Assumption \ref{assump2}.
\begin{lemma}
Under Assumption \ref{assumption_freq}, the mapping \eqref{integral1}satisfies Assumption \ref{assump2} with $\varphi:=[x^\top,\tau]^\top$, ${\chi}_{bl}:=\mu$, $\delta:=a$, $f^{\delta}_{\varphi}:=F_x\times F_{\tau}$, and  $f_{\chi}=R\mu$. 
\end{lemma}
\textbf{Proof:} Using the definitions of the lemma, and the result of Lemma \ref{lemma_dither}, for each $a>0$ there exists a $\tilde{T}>0$ such that for each $(x,\tau)\in \left(\mathbb{R}^n\times\mathcal{T}_{C}\right) \cap K$ with $K\subset\mathbb{R}^{n+1}$ compact, the following holds for each case:
\begin{equation}
\frac{1}{N\tilde{T}}\int_{0}^{N\tilde{T}} \left[\left(\begin{array}{c}
F_{x}(x,\tau,\mu(s))\\
F_{\tau}
\end{array}\right)-f^{a}_A(x,\tau)\right]ds=0,
\end{equation}
for all $N\in\mathbb{Z}_{\geq0}$. Since any $L\in\mathbb{R}_{>0}$ can be written as $L=N\tilde{T}+\tilde{L}$ where $|\tilde{L}|\leq \tilde{T}$, it suffices to consider the integral over the interval $[0,\tilde{L}]$.
%
%
%
By the proof of Lemma \ref{lemma_dither} in the Appendix, the integrals of $\mu(s)$ and $\mu(s)\mu(s)^\top$ are bounded on any finite time. Thus, by the construction of $F_x$ and the fact that $(x,\tau)\in K$ with $K$ compact, and the continuity of $f^{a}_A(x,\tau)$ and $H_j$, as well as the smoothness of $f$, there exists $M_{K,a}>0$ such that 
\begin{equation}
\left|\int_{0}^{\tilde{L}} \left[\left(\begin{array}{c}
F_{x}(x,\tau,\mu(s))\\
F_{\tau}
\end{array}\right)-f^{a}_A(x,\tau)\right]ds\right|\leq M_{K,a}.
\end{equation}
The bound \eqref{bound_average} holds with $L=N\tilde{T}+\tilde{L}$ and $\sigma_{K,\delta}(L)=M_{k,a}/L$.  \null \hfill \null $\blacksquare$
%
%
\subsection{Stability Analysis of Average Systems}
\label{subsection_stability_analysis}
Having obtained a well-defined average mapping $f_A^{a}$, we now define the average jump map as
\begin{equation}\label{jumps}
G_A(y):=G_x(y)\times\{T_{\min}\},
\end{equation}
with $G_x$ defined as in Sections \ref{subsection:1}-\ref{subsection:4}. This definition is consistent with Definition \ref{definition_average_system}. The average hybrid dynamics of the HDS $\mathcal{H}_{a,\varepsilon}$ are defined as
\begin{subequations}\label{average_system0}
\begin{align}
y\in C&:=\mathbb{R}^{n+m}\times \mathcal{T}_C,~~~~~~~~~\dot{y}=f_A^a(y)\label{average1a}\\
y\in D&:=\mathbb{R}^{n+m}\times  \mathcal{T}_D,~~~~~~~y^+=G_A(y).\label{average1b}
\end{align}
\end{subequations}
We now show that for each Case 1-4, system \eqref{average_system0} satisfies Assumption \ref{assumption_stability}.

\textbf{Case 1:} The average system has state $y\in\mathbb{R}^{2n+1}$ and dynamics \eqref{average_system0}, where $f_A^a$ is defined as in \eqref{integral1}, $G_A$ is defined as in \eqref{jumps}, and $\mathcal{T}_C,\mathcal{T}_D$ are defined as in \eqref{sets_1}. For this HDS, the following Lemma holds.
\begin{lemma}\label{lemma_1}
Let $\mathcal{A}_x$ be given by \eqref{optimal_1}. Under Assumption \ref{assumption1}, the HDS \eqref{average_system0} renders the set $\mathcal{A}_{\varphi}:=\mathcal{A}_x\times\mathcal{T}_C$  SGPAS as $a\to0^+$.
\end{lemma}

\textbf{Proof:} In the first step of the proof we neglect the $\mathcal{O}(a)$ perturbation term in the mapping \eqref{integral1}, and we establish UGAS of the set $\mathcal{A}_x\times\mathcal{T}_C$. In the second step, we use the robustness properties of well-posed hybrid systems (e.g., \cite[Thm. 7.21]{Goebel:12}) to establish SGPAS as $a\to0^+$ for the original $O(a)$-perturbed system. 

\textsl{Step 1:} First, we consider the case when $k_1=0$ and $F_{\tau}=0.5$. Consider the Lyapunov function
\begin{equation}\label{lyapunov1}
 V_{k_1}(y)=\frac{1}{4}|y_2-y_1|^2+\frac{1}{4}|y_2|_{\mathcal{A}_{\phi}}^2+k_2y_3^2(\phi(y_1)-\phi^*), 
\end{equation}
which is radially unbounded and positive definite with respect to $\mathcal{A}$. Using the definition of $\dot{y}$, and denoting as $z^*=P_{\mathcal{A}_f}(y_2)$ the Euclidean projection of $y_2$ on $\mathcal{A}_{\phi}$, we obtain that $\dot{V}_{k_1}(y)=\nabla V_{k_1}^\top \dot{y}$ satisfies
\begin{align}
\dot{V}_{k_1}\leq-&\frac{1}{\tau}|y_2-y_1|^2\notag\\
&-k_2\tau\left(\nabla \phi(y_1)^\top(y_1-z^*)-\phi(y_1)+\phi^*\right),\label{decrease_derivative_flows}
\end{align}
which, by convexity, implies that $\dot{V}_{k_1}\leq0$ for all $y\in C$. Moreover, by Lemma A.1 in the Appendix, when $\mathcal{A}_{\phi}$ is a singleton we have that $\dot{V}_{k_1}<0$ for all $y\in C\backslash\mathcal{A}$. On the other hand, when $\mathcal{A}_{\phi}$ is not a singleton, but $\nabla \phi$ is globally Lipschitz, the right hand side of \eqref{decrease_derivative_flows} can be further upper bounded as $\dot{V}_{k_1}\leq-\frac{1}{\tau}|y_2-y_1|^2-\frac{k_2\tau}{2L}|\nabla \phi(y_1)|^2$, and by the definition of $\mathcal{A}$ and convexity, this implies that $\dot{V}_{k_1}<0$ for all $y\in C\backslash\mathcal{A}$. In addition, during jumps the change of the Lyapunov function $\Delta V_{k_1}(y):=V_{k_1}(y^+)-V_{k_1}(y)$  satisfies
\begin{align}\label{non_increase_jumps}
\Delta V_{k_1}(y)=-k_2(\phi(y_1)-\phi(z^*))(y_3^2-T^2_{\min})\leq 0,
\end{align}
for all $y\in D$. Inequalities \eqref{decrease_derivative_flows} and \eqref{non_increase_jumps} imply that $\mathcal{A}$ is stable. Since the Lyapunov function does not increase during jumps, and it is strictly decreasing during flows, it follows that for all $c>0$ there is no complete solution $y$ such that $V_{k_1}(y(t,j))=c$ for all $(t,j)\in\text{dom}(\phi)$. Therefore, by the Hybrid Invariance Principle \cite[Thm. 8.8]{Goebel:12} there exists a $\beta_1\in\mathcal{K}\mathcal{L}$ such that the set $\mathcal{A}_{\varphi}$ is UGpAS.

Next, let $k_1\geq0$ and $F_{\tau}=1$. Consider the candidate Lyapunov function
\begin{equation}\label{LyapunovLemma92}
V_{k_1}(y)=\frac{1}{2}|y_2-z^*|^2+k_2y_3^2(\phi(y_1)-\phi^*),
\end{equation}
which is radially unbounded and positive definite with respect to the set $\mathcal{A}_x$. The time derivative satisfies $\dot{V}_{k_1}=-2\tau k_2(\phi^*-\phi(y_1)+\nabla \phi(y_1)^\top (y_1-z^*))-k_2k_1\tau^2|\nabla\phi(y_1)|^2\leq 0$, where the last inequality follows by convexity. During jumps the Lyapunov function $V_{k_1}$ also satisfies inequality \eqref{non_increase_jumps}. Since by assumption the set $\mathcal{A}_x$ is a singleton, by Lemma A.1 in the Appendix, and the hybrid invariance principle \cite[Thm. 8.8]{Goebel:12} every maximal solution converges to the largest weakly invariant subset contained in the set $W_r=\{y\in\mathbb{R}^{2n+1}:y_1=z^*,|y_2-z^*|=\sqrt{2r},y_3\in\mathcal{T}_C\cup\mathcal{T}_D\}$. Since during flows any solution in $W_r$ satisfies $y_1=z^*\implies \dot{y}_1=0=2y_3^{-1}(y_2-z^*)$ and $\dot{y}_2=0$, it must be the case that $y_2=z^*$. Since $y_1^+=y_1$ and $y_2^+=y_2$, the largest weakly invariant subset of $W_r$ is $W_0=\mathcal{A}_x$. Therefore, by the Hybrid Invariance Principle, there exists a $\beta_1\in\mathcal{K}\mathcal{L}$ such that the set $\mathcal{A}_{\varphi}$ is UGpAS. 

\textsl{Step 2:} Since the HDS \eqref{average_system0} with $O(a)=0$ satisfies the Basic Conditions, and $\mathcal{A}_{\varphi}$ is UGAS, by \cite[Thm. 7.21]{Goebel:12} the original $O(a)$-perturbed system \eqref{average_system0} renders the set $\mathcal{A}_{\varphi}$ SGPpAS as $a\to 0^+$ with $\beta_1\in\mathcal{K}\mathcal{L}$. Since such solutions satisfy the bound \eqref{wdelta},  this rules out finite escape times, and since $G_{es}(D_{es})\subset C_{es}\backslash D_{es}$,  by Lemma \ref{existence_completness} the solutions are complete. This establishes the result of the lemma.
 \null \hfill \null $\blacksquare$

\textbf{Case 2:} The average system has state $y\in\mathbb{R}^{2n+1}$ and dynamics \eqref{average_system0}, where $f_A^a$ is defined as in \eqref{integral1} with $k_1=0$, $G_A$ is defined as in \eqref{jumps}, and $\mathcal{T}_C,\mathcal{T}_D$ are defined as in \eqref{sets_2}. For this HDS, the following Lemma holds.
\begin{lemma}\label{lemma_2}
Let $\mathcal{A}_{x}$ be given by \eqref{optimal_1}. Under Assumption \ref{assumption2}, system \eqref{average_system0} renders the set $\mathcal{A}_{\varphi}:=\mathcal{A}_x\times\mathcal{T}_C$ SGPES as $a\to0^+$. 
\end{lemma}
\textbf{Proof:}  We follow the same two steps as in the proof of Lemma \ref{lemma_1}, and we use the fact that $|y|^2_{\mathcal{A}_{\varphi}}=|y_1-z^*|^2+|y_2-z^*|^2$ for all $y\in C\cup D$.

\textsl{Step 1:} Setting $k_1=0$, neglecting the $O(a)$ perturbation in \eqref{integral1}, using the Lyapunov function \eqref{lyapunov1} (we omit the sub-index $k_1$), and using strong convexity and globally Lipschitz of $\nabla \phi$, we obtain
\begin{equation}\label{useful_bounds4}
\underline{c}|y|^2_{\mathcal{A}_{\varphi}}\leq V(y)\leq \overline{c}|y|^2_{\mathcal{A}_{\varphi}},
\end{equation}
%
with $\underline{c}:=0.25\min\{1,2kT_{\min}^2\theta\}$ and $\bar{c}:=\max\{3,6kT_{\max}^2L\}$. 
%
%
%
Since the flow map is the same as in Case 1 with $k_1=0$, the derivative $\dot{V}$ satisfies \eqref{decrease_derivative_flows}. Thus, using strong convexity and the definition of $\mathcal{T}_C$ we obtain $\dot{V}(y)\leq -\rho|y|^2_{\mathcal{A}_{\varphi}}$, where $\rho:=\min\{0.5\tilde{\rho},0.25kT_{\min}\theta\}$. Using the upper bound of \eqref{useful_bounds4}, defining $\gamma:=1-\frac{T^2_{\min}}{T_{\max}^2}-\frac{1}{2k\theta T_{\max}^2}$, which due to \eqref{dwell_time} satisfies $\gamma\in(0,1)$, and $\lambda_2:=\min\left\{\frac{\rho}{\bar{c}},-\log(1-\gamma)\right\}$, we finally obtain:
\begin{align}
\dot{V}(y)&\leq-\lambda_2 V(y),~~~\forall~y\in C.~\label{strict_flow3s}
\end{align}
During jumps, the Lyapunov function satisfies 
%
%
%
\begin{align*}
\Delta V(z^+)\leq&-\frac{1}{4}|y_2-y_1|^2-\frac{1}{4}|y_2-z^*|^2\notag\\
&-k\tau^2(\phi(y_1)-\phi^*)\left(1-\frac{T_{\min}^2}{\tau^2}-\frac{1}{2k\theta\tau^2} \right),
\end{align*}
where we used the strong convexity of $\phi$. By the definition of $\gamma$, the fact that $0<\gamma<1$, and that $\tau=T_{\max}$ in the jump set, we obtain $\Delta V(z)\leq-\gamma V(y)$, for all $y\in D$. Thus, during jumps the Lyapunov function satisfies
\begin{align}\label{jump_final12}
V(y^+)\leq (1-\gamma)~V(y)\leq \exp(-\lambda_2)~V(y),
\end{align}
for all $y\in D$. Inequalities \eqref{useful_bounds4}, \eqref{strict_flow3s}, and \eqref{jump_final12}, imply that the HDS \eqref{average_system0} with $O(a)=0$ renders the set $\mathcal{A}_{\varphi}$ UGpES. Using again \cite[Thm. 7.21]{Goebel:12} we obtain SGPpES as $a\to0^+$ for the $O(a)$-perturbed system  \eqref{average_system0}. The $\mathcal{K}\mathcal{L}$ bound \eqref{wdelta} rules out finite escape times, and since the jumps always happen when $\tau=T_{\max}$, and always reset the timer to $T_{\min}$, maximal solutions do not stop. Therefore, by Lemma \ref{existence_completness} every maximal solution satisfying \eqref{wdelta} is complete. \null \hfill \null $\blacksquare$

\textbf{Case 3:} The average system has a state $y\in\mathbb{R}^{n+m+1}$, and hybrid dynamics \eqref{average_system0}, where $f_A^a$ is defined as in \eqref{integral1}, $G_A$ is defined as in \eqref{jumps}, and $\mathcal{T}_C,\mathcal{T}_D$ are defined as in \eqref{sets_3}. For this HDS we have the following lemma.

\begin{lemma}\label{lemma_3}
Let $\mathcal{A}_x$ be given by \eqref{optimal_lag1} with Lagrangian \eqref{lagrangian1}. Under Assumption \ref{assumption3}, the HDS system \eqref{average_system0} renders the set $\mathcal{A}_{\varphi}:=\mathcal{A}_x\times\{T_{\min}\}$ SGPES as $a\to0^+$.
\end{lemma}

\textbf{Proof:} Since $y_3(t)=T_{\min}$ for all $t\geq0$, and the dynamics of $(y_1,y_2)$ and $y_3$ are uncoupled, it suffices to study the properties of $y_1,y_2$ with respect to $\mathcal{A}_x$. Neglecting the $O(a)$-perturbation, and ignoring the jumps, by \cite[Lemma 1]{GuananLetters} there exists a Lyapunov function of the form  $V(\tilde{y})=\tilde{y}^\top P_3\tilde{y},~~P_3>0,~~\tilde{y}:=[(y_1-z^*)^\top,(y_2-\lambda^*)^\top]^\top$. By  \cite[Lemma 2]{GuananLetters} there exists $\lambda_{3} \in\mathbb{R}_{>0}$ such that  $\dot{V}(\tilde{y})\leq -\lambda_{3} V(\tilde{y})$. Therefore, the point $y=z^*$ is UGpES. Since the average HDS with $O(a)=0$ is well-posed, \cite[Thm. 7.21]{Goebel:12} establishes SGPpES as $a\to0^+$ for the original $O(a)$-perturbed average system. Since the system has no jumps and the flow set is unbounded, every maximal solution satisfying \eqref{wdelta} is complete.   \null \hfill \null $\blacksquare$

\textbf{Case 4:} Let $\mathcal{A}_x$ be given by \eqref{optimal_lag1} with Lagrangian \eqref{lagrangian2}. The average system has state $y\in\mathbb{R}^{n+m+1}$, and hybrid dynamics \eqref{average_system0}, where $f_A^a$ is defined as in \eqref{integral1}, $G_A$ is defined as in \eqref{jumps}, and $\mathcal{T}_C,\mathcal{T}_D$ are defined as in \eqref{sets_3}. For this HDS we have the following lemma.
\begin{lemma}\label{lemma_4}
Under Assumption \ref{assumption3}, the system \eqref{average_system0} renders the set $\mathcal{A}_{\varphi}:=\mathcal{A}_x\times\{T_{\min}\}$ SGPES as $a\to0^+$.
\end{lemma}
\textbf{Proof:}  Since $y_3(t)=T_{\min}$ for all $t\geq0$ and the dynamics of $(y_1,y_2)$ and $y_3$ are uncoupled, it suffices to study the properties of $y_1,y_2$ with respect to $\mathcal{A}_x$. Neglecting the $O(a)$-perturbation, by \cite[Lemmas 3 $\&$ 4]{GuananLetters} there exists a quadratic Lyapunov function  $V=\tilde{y}^\top P_4\tilde{y},~~P_4>0,~~\tilde{y}=[(y_1-z^*)^\top,(y_2-\lambda^*)^\top]^\top$ and some $\lambda_4\in\mathbb{R}_{>0}$ such that  $\dot{V}(y)\leq -\lambda_{4}V(y)$, during flows of the system. Therefore, the point $y=z^*$ is UGpES. Since the average HDS with $O(a)=0$ is well-posed, \cite[Thm. 7.21]{Goebel:12} establishes SGPpES as $a\to0^+$ for the original $O(a)$-perturbed average system. Since the system has no jumps and the flow set is unbounded, every maximal solution satisfying \eqref{wdelta} is complete.   \null \hfill \null $\blacksquare$.

Lemmas \ref{lemma_1}, \ref{lemma_2}, \ref{lemma_3}, and \ref{lemma_4} imply that the average HDS \eqref{average_system0} satisfies Assumption \ref{assumption_stability} for Cases 1-4. Thus, by Theorem \ref{main_theorem_averaging}, the original hybrid dynamics \eqref{HAES} render the set $\mathcal{A}:=\mathcal{A}_x\times \mathcal{T}_C\times \mathbb{T}^n$ SGPpAS as $(a,\varepsilon)\to0^+$ with the same $\mathcal{K}\mathcal{L}$ function as the average system. Completeness of solutions follows by the absence of finite escape times implied by the bound \eqref{wdelta}, the fact that no solution can leave $\mathcal{T}_C\cup\mathcal{T}_D$ due to jumps or flows, and Lemma \ref{existence_completness} in the Appendix.  Finally, since all the HAES considered in Theorems \ref{thm1a}-\ref{thm4} satisfy the Basic Conditions and render the compact set $\mathcal{A}$ SGPAS as $(a,\varepsilon)\to0$, Structural Robustness follows now directly by Proposition \ref{robustness_proposition} in the Appendix. This establishes the stability, robustness, and completness results of Theorems \ref{thm1a}, \ref{thm2}, \ref{thm3}, and \ref{thm4}.
\subsection{Convergence Bounds}
We now proceed to establish the convergence bounds for Cases 1-4. To shorten notation we use $\tilde{\phi}=\phi-\phi^*$ to denote the sub-optimality measure.

\noindent
\textbf{Case 1. Step 0:}
 Let $K_0$ and $\nu$ be given. Without loss of generality we assume that $\nu\in(0,1)$. Consider the average nominal hybrid system $\mathcal{H}^A$ \eqref{average_system0} with state $y$ and $O(a)$-perturbation set to zero. Consider also the average $O(a)$-perturbed hybrid system $\mathcal{H}_a^A$ given by \eqref{average_system0}  with state $y^a$, as well as original hybrid system $\mathcal{H}$ with state $\tilde{x}=[x^\top,\mu^\top,\tau]^\top$. Define the set $\tilde{K}_0=K_0+\mathbb{B}$ and the quantity $m:=\max_{x_0\in \tilde{K}_0}|x_0|_{\mathcal{A}_x}$. Let $\beta_1$ be the $\mathcal{K}\mathcal{L}$-bound that characterizes the UGAS property of $\mathcal{H}^A$, which was established in Step 1 of the proof of Lemma \ref{lemma_1}. Define the set $K_1=\left\{x\in\mathbb{R}^{2n}:|x|_{\mathcal{A}_x}\leq \beta_1 \left(m,0\right)+1\right\}$. Since $K_1$ is compact, there exists $M>0$ such that $K_1\subset M\mathbb{B}$. Define $M_1:=\max_{z\in M\mathbb{B}} |\phi(z)-\phi^*|$ and \tcbb{$\hat{\lambda}:=\frac{5}{T_{\min}^2}\max\{k_2^{-1},T_{\max}^2,T_{\max}M_1, M_1,1\}$}. Let $\nu'\in(0,0.5\nu\min\{1/\hat{\lambda},1\})$.  By continuity of $\phi$ and the functions $V_{k_1}$ in \eqref{lyapunov1} and \eqref{LyapunovLemma92}, there exists $\delta_1^*\in(0,\min\{\nu',2(T_{med}-T_{\min})\})$ such that $|r_1-r_2|\leq \frac{\delta_1}{4}\implies |\phi(r_1)-\phi(r_2)|\leq \frac{\nu'}{4}$ and $|V_{k_1}(r_1)-V_{k_1}(r_2)|\leq \frac{\nu'}{4}$, for all $r_1,r_2\in M\mathbb{B}$ and for all $\delta_1\in(0,\delta_1^*)$. Fix $\delta_1\in(0,\delta_1^*)$ and let $T^*>1$ be such that for all $w\geq T^*$ we have $\beta(m,w)<\frac{\delta_1}{4}$.  Let $\overline{T}:=T^*+1$. Then, by their respective properties of UGAS, SGPAS as $a\to0^+$, and SGPAS as $(a,\varepsilon)\to0^+$, and by Proposition \ref{closeness_on_compact_time_domain}, there exists $a^*\in(0,\frac{\delta_1}{4})$ such that for all $a\in(0,a^*)$ there exists $\varepsilon^*\in(0,1)$ such that for all $\varepsilon\in(0,\varepsilon^*)$ the following properties hold: \textbf{(a)} All solutions $y$ of $\mathcal{H}^A$ with $(y_1(0,0),y_2(0,0))\in \tilde{K}_0$ satisfy $y(\tilde{t},j)\in M\mathbb{B}\times\mathcal{T}_C$ for all $(\tilde{t},j)\in\text{dom}(y)$. \textbf{(b)} All solutions $y_a$ of $\mathcal{H}_a^A$ with $(y_{a,1}(0,0),y_{a,2}(0,0))\in \tilde{K}_0$ satisfy $y_a(\hat{t},j)\in M\mathbb{B}\times\mathcal{T}_C$ for all $(\hat{t},j)\in\text{dom}(y_a)$. \textbf{(c)} All solutions $\tilde{x}$ of $\mathcal{H}$ with $x(0,0)\in \tilde{K}_0$ satisfy $\tilde{x}(t,j)\in M\mathbb{B}\times\mathcal{T}_C\times\mathbb{T}^n$ for all $(t,j)\in\text{dom}(\tilde{x})$, and $|x(t,j)|_{\mathcal{A}_x}\leq \frac{\delta_1}{4}$ for all $t+j\geq T^*$.\textbf{(d)} For each solution $y_a$ of $\mathcal{H}^A_a$ with $y_a(0,0)\in K_0+a\mathbb{B}$ there exists a solution $y$ of $\mathcal{H}^A$ with $y(0,0)\in K_0$ that is $(\overline{T},\frac{\delta_1}{4})$-close. \textbf{(e)} For each solution $\tilde{x}$ of $\mathcal{H}$ with $x(0,0)\in K_0$ there exists a solution of $\mathcal{H}^A_a$ with $y_a\in K_0$ that is $(\overline{T},\frac{\delta_1}{4})$ close to $(x,\tau)$. 
 
By the second part of Property \textbf{(c)}, the continuity of $\phi$ and $V_{k_1}$ on $M\mathbb{B}$, and the choice of $\delta_1$, we obtain 
\begin{equation}\label{steady_state_error}
\tilde{\phi}(x_1(t,j))\leq \frac{\nu'}{4},~~\forall~(t,j)\in\text{dom}(\tilde{x}),~t+j\geq T^*,
\end{equation}
and $\lim\sup_{t+j\to\infty} V(\tilde{x}(t,j))\leq \nu$.

\textbf{Step 1:} From the stability analysis of Step 1 in the proof of Lemma \ref{lemma_1}, we known that $V_{k_1}(y(\tilde{t},j))\leq V_{k_1}(y(\tilde{t}',j))$ for all $\tilde{t}\geq \tilde{t}'$ such that $(\tilde{t},j),(\tilde{t}',j)\in\text{dom}(y)$. Let $\tilde{s}_j=\min\{\tilde{t}: (\tilde{t},j)\in\text{dom}(y)\}$. Using the structure of the Lyapunov functions \eqref{lyapunov1} and \eqref{LyapunovLemma92}, for each $j$ such that $(\tilde{t},j)\in\text{dom}(y)$ the following bound holds:
\begin{equation}\label{inequaLyapunov}
\tcbb{k_2y_3^2(\tilde{t},j)\tilde{\phi}(y_1(\tilde{t},j))\leq V_{k_1}(y(\tilde{s}_j,j)),~~\forall~\tilde{t}>\tilde{s}_j.} 
\end{equation}
By Property (\textbf{d}), for each solution of $\mathcal{H}_a^A$ with $y_a(0,0)\in K_0+a\mathbb{B}$ there exists a solution $y$ of $\mathcal{H}^A$ with $y(0,0)\in K_0$ that satisfies the following: For all $(\hat{t},j)\in\text{dom}(y_a)$ with $\hat{t}+j<\overline{T}$, there exists $\tilde{t}'$ such that $(\tilde{t}',j)\in\text{dom}(y)$, $|\hat{t}-\tilde{t}'|\leq \frac{\delta_1}{4}$ and $|y_a(\hat{t},j)-y(\tilde{t}',j)|\leq \frac{\delta_1}{4}$. Using again the uniform boundedness of $y$ and $y_a$, the continuity of $\phi$ on $M\mathbb{B}$, and the choice of $\delta_1$, we have that $|y_{1,a}(\hat{t},j)-y_1(\tilde{t}',j)|\leq \frac{\delta_1}{4}~\implies \phi(y_{1,a}(\hat{t},j))\leq \phi(y_1(\tilde{t}',j))+\frac{\nu}{4}$ and $y_{3,a}(\hat{t},j)-\frac{\nu'}{4}\leq y_3(\tilde{t}',j)$. \tcbb{The last inequality also implies $y^2_{3,a}(\hat{t},j)\leq y_3^2(\tilde{t}',j)+\nu' N_1$, where $N_1:=2T_{\max}\frac{1}{4}+\frac{1}{16}$.} Therefore, the left hand side of \eqref{inequaLyapunov} evaluated at $\tilde{t}'$ can be lower bounded as:
\begin{align}\label{lower_bound_proofTransient1}
&\tcbb{k_2y_{3,a}^2(\hat{t},j)\tilde{\phi}(y_{1,a}(\hat{t},j))\leq}\tcbb{k_2y_3^2(\tilde{t}',j)\tilde{\phi}(y_1(\tilde{t}',j))+k_2\nu' N_2,}
\end{align}
\tcbb{where $N_2:=\frac{1}{4}T_{\max}^2+N_1M_1+\frac{1}{4}$.}  Since Property $\textbf{(d)}$ also implies $|\hat{s}_j-\tilde{s}_j|\leq \frac{\delta_1}{4}$ and $|y(\tilde{s}_j,j)-y_a(\hat{s}_j,j)|\leq \frac{\delta_1}{4}$, where $\hat{s}_j=\min\{\hat{t}: (\hat{t},j)\in\text{dom}(y_a)\}$, using the continuity of $V_{k_1}$ we obtain $V_{k_1}(y(\tilde{s}_j,j)) \leq V_{k_1}(y_{a}(\hat{s}_j,j))+\frac{\nu'}{4}$. Using this inequality to upper bound \eqref{inequaLyapunov}, and using \eqref{lower_bound_proofTransient1}:
\begin{equation}\label{final_bound_step1}
\tcbb{k_2y_{3,a}^2(\hat{t},j)\tilde{\phi}(y_{1,a}(\hat{t},j))\leq V_{k_1}(y_{a}(\hat{s}_j,j))+\nu' N_3,}
\end{equation}
\tcbb{where $N_3=\frac{1}{4}+k_2N_2$,} which holds for all $\hat{t}>\hat{s}_j$ such that $(\hat{t},j)\in\text{dom}(y_a)$ and $\hat{t}+j\leq \overline{T}$. 

\noindent
\textbf{Step 2:} We now repeat the exact same procedure of Step 1 using Property (\textbf{e}) instead of (\textbf{d}) to relate the properties of $y_a$ with the properties of $\tilde{x}$. Indeed, by Property (\textbf{e}), for each solution $\tilde{x}$ of $\mathcal{H}$ with $x(0,0)\in K_0$ there exists a solution $y_a$ of $\mathcal{H}_a^A$ with $(y_{1,a}(0,0),y_{2,a}(0,0))\in K_0$ that satisfies the following: For all $(t,j)\in\text{dom}(\tilde{x})$ with $t<\overline{T}$, there exists $\hat{t}'$ such that $(\hat{t}',j)\in\text{dom}(y_a)$, $|t-\hat{t}'|\leq \frac{\delta_1}{4}$ and $|\tilde{x}(t,j)-\tilde{y}_a(\hat{t}',j)|\leq \frac{\delta_1}{4}$. Using this property, and similarly to Step 1, we can lower bound the left-hand side of \eqref{final_bound_step1} evaluated at $\hat{t}'$ as follows:
\begin{align}\label{usefulboudstep2}
&\tcbb{k_2\tau^2(t,j)\tilde{\phi}(x_{1}(t,j))\leq}\tcbb{k_2y_{3,a}^2(\hat{t}',j)\tilde{\phi}(y_{1,a}(\hat{t}',j))+k_2\nu' N_2},
\end{align}
and we can upper bound the right hand side of \eqref{final_bound_step1} as $V_{k_1}(y_{a}(\hat{s}_j,j))+\nu' N_3\leq V_{k_1}(\tilde{x}(\underline{t}_j,j))+\frac{\nu'}{4}+\nu' N_3$, where $\underline{t}_j=\min\{t: (t,j)\in\text{dom}(\tilde{x})\}$. Using this bound, as well as \eqref{final_bound_step1} and \eqref{usefulboudstep2}, we obtain:
 \begin{equation}\label{almostthere2}
\tcbb{k_2\tau^2(t,j)\tilde{\phi}(x_{1}(t,j))\leq V_{k_1}(\tilde{x}(\underline{t}_j,j))+\nu'N_4,}
\end{equation}
\tcbb{where $N_4=N_3+k_2N_2+\frac{1}{4}$, and which holds for all $t>\underline{t}_j$, such that $(t,j)\in\text{dom}(\tilde{x})$ and $t+j\leq \overline{T}$.  Dividing both sides by $k_2\tau^2$, and using $\tau^2\geq T_{\min}$, we obtain:}
\begin{equation}\label{almostthere3}
\tcbb{\tilde{\phi}(x_{1}(t,j))\leq \frac{V_{k_1}(\tilde{x}(\underline{t}_j,j))}{k_2\tau^2(t,j)}+\nu' \frac{N_4}{k_2T_{\min}^2}.}
\end{equation}
Since $\dot{\tau}\in\{0.5,1\}$, and $\tau(\underline{t}_{j},j)=T_{\min}$ for all $j$, it follows that $\tau^2(t,j)\geq(0.5(t-\underline{t}_j)+T_{\min})^2> 0.25(t-\underline{t}_j)^2$. Substituting in \eqref{almostthere3} we obtain:
\begin{equation}\label{almosttherefinal}
\tcbb{\tilde{\phi}(x_{1}(t,j))\leq \frac{4V_{k_1}(\tilde{x}(\underline{t}_j,j))}{k_2(t-\underline{t}_j)^2}+\nu' \frac{N_4}{k_2T_{\min}^2},}
\end{equation}
which holds for all $t>\underline{t}_j$, such that $(t,j)\in\text{dom}(\tilde{x})$ and $t+j\leq \overline{T}$.  Since $a<\frac{\delta_1}{4}$, $|\mu(t,j)|\leq 1$ for all $(t,j)\in\text{dom}(\tilde{x})$, by continuity of $\phi$ we obtain $\phi(x_1+a\mu)\leq \phi(x_1)+\frac{\nu}{4}$ for all $|x_1|\leq M$. Combining this inequality with \eqref{almosttherefinal} and \eqref{steady_state_error}, and using the facts that $N_4/k_2T_{\min}^2\leq\hat{\lambda}$ and $\nu'\hat{\lambda}<0.5\nu$, we obtain the bound \eqref{inequality_algo1}, which holds for all $t>\underline{t}_j$ such that $(t,j)\in\text{dom}(\tilde{x})$.   \null \hfill \null $\blacksquare$

\textbf{Cases 2-4: Exponential Decrease.}
The SGPES results of Lemmas \ref{lemma_2}, \ref{lemma_3} and \ref{lemma_4},  imply that in each Case $i\in\{2,3,4\}$, every solution of the average HDS satisfies the bound \eqref{KL_average_system} with $\omega(y)=|y|_{\mathcal{A}_{\varphi}}$, and $\beta_i\in\mathcal{K}\mathcal{L}$ given by $\beta_i(r,s)=\alpha_{1,i}\exp \left(-\alpha_{2,i} s\right)r$, where $\alpha_{1,i},\alpha_{2,i}>0$, $i\in\{2,3,4\}$, which establishes the exponential convergence result for Cases 2-4.  

To establish inequality \eqref{exponential_decrease_during_jumps} for Case 2, we follow similar steps as in the proof of the convergence bound of Case 1. In particular, by the proof of Lemma \ref{lemma_2}, we have that $\dot{V}(y(t,j))\leq 0$ during flows, where $V$ is given by \eqref{lyapunov1}, which also satisfies inequality \eqref{inequaLyapunov}.  Let $\tilde{s}_j:=\min\{\tilde{t}:(\tilde{t},j)\in\text{dom}(y)\}$. Since the initialization and the jump rule imply that $y_2(\tilde{s}_j,j)=y_1(\tilde{s}_j,j)$ and $y_3(\tilde{s}_j,j)=T_{\min}$, for all $j\in\mathbb{Z}_{\geq0}$, using strong convexity to upper bound the right-hand side of \eqref{inequaLyapunov} we obtain:
\begin{align}\label{first_forrecursion}
\tcbb{\tilde{\phi}(y_1(\tilde{t},j))\leq \tilde{\alpha}_0 \tilde{\phi}(y_1(\tilde{s}_j,j)),}
\end{align}
for all $(\tilde{t},j)\in\text{dom}(y)$ such that $\tilde{t}\geq \tilde{s}_j$, where $\tilde{\alpha}_0=\alpha_0\tilde{\gamma}$. Now, since the jump rule does not change $y_1$, it follows that $\tilde{\phi}(y_1(\tilde{s}_{j+1},{j+1}))=\tilde{\phi}(y_1(\tilde{s}_{j+1},{j}))$. Moreover, since $\tilde{s}_{j+1}=\tilde{s}_{j}+\Delta T$, and since at the end of the periods of flow we have $y_3(\tilde{s}_{j}+\Delta T,j)=T_{\max}$, it follows that
\begin{align}\label{second_forrecursion}
\tcbb{\tilde{\phi}(y_1(\tilde{s}_{j+1},{j+1}))\leq \tilde{\gamma} \tilde{\phi}(y_1(\tilde{s}_j,j)),}
\end{align}
which holds for all $j\in\mathbb{Z}_{\geq0}$. Combining \eqref{first_forrecursion} and \eqref{second_forrecursion}, we can now
follow the same steps as in Case 1 using closeness of solutions between the average system and the HAES. In particular, we use the same construction of Step 0 in Case 1, but now we use $\nu'\in(0,(1-\tilde{\gamma}/\tilde{\alpha}_0)\nu)$, and we choose again $a$ and $\varepsilon$ such that all properties \textbf{(a)}-\textbf{(e)} hold, including the bound \eqref{steady_state_error}. Then, by using again closeness of solutions between the trajectories $y_a$ of the perturbed average hybrid system $\mathcal{H}_a^A$, the trajectories $y$ of the unperturbed  average hybrid system $\mathcal{H}^A$, and the trajectories of the HAES we obtain that every solution of the HAES satisfies the bound
\begin{equation}\label{transient_bound}
\tcbb{\tilde{\phi}(x_{1}(t,j))\leq \tilde{\alpha}_0 \tilde{\phi}(x_{1}(\underline{t}_j,j))+\frac{\nu'}{4},}
\end{equation}
for all $t>\underbar{t}_j$ such that $t+j\leq \overline{T}$, and
\begin{equation}\label{inequlityforiterations}
\tcbb{\tilde{\phi}(x_{1}(\underline{t}_{j+1},j+1))\leq \tilde{\gamma} \tilde{\phi}(x_{1}(\underline{t}_j,j))+\frac{\nu'}{4},}
\end{equation}
for all $j\geq0$ such that $t+j\leq \overline{T}$. \tcbb{We can now use inequality \eqref{inequlityforiterations} to iterate over $j$ starting with $j=0$:}
\begin{align*}
\tcbb{\tilde{\phi}(x_{1}(\underline{t}_{j},j))}&\tcbb{\leq\tilde{\gamma}^j \tilde{\phi}(x_1(0,0))+\frac{\nu'}{4}\frac{1-\tilde{\gamma}^j}{1-\tilde{\gamma}}.}
\end{align*}
\tcbb{Using this expression to upper bound \eqref{transient_bound}, and the facts that $\tilde{\gamma}\in(0,1)$ and $\nu'\tilde{\alpha}_0< (1-\tilde{\gamma})\nu$, we finally obtain:}
\begin{equation}\label{last_one_exponential}
\tcbb{\tilde{\phi}(x_{1}(t,j))\leq \tilde{\alpha}_0\tilde{\gamma}^j\tilde{\phi}(x_1(0,0))+\frac{\nu}{2},}
\end{equation}
which holds for all $(t,j)\in\text{dom}(\tilde{x})$ such that $t>\underline{t}_j$ and $t+j\leq \overline{T}$. Since $a<\frac{\delta_1}{4}$, $|\mu(t,j)|\leq 1$, and using the uniform continuity of $\phi$ on $M\mathbb{B}$, we obtain $\phi(x_1+a\mu)\leq \phi(x_1)+\frac{\nu}{4}$. \tcbb{Combining this inequality with \eqref{last_one_exponential} and \eqref{steady_state_error}, and using $\tilde{\alpha}_0<\alpha_0$, we obtain the desired bound \eqref{exponential_decrease_during_jumps}.}\null \hfill \null $\blacksquare$

\section{Conclusions and Outlook}
\label{sec_conclusions}
This paper presents a new class of zero-order optimization dynamics with acceleration and restarting mechanisms. These algorithms can be modeled as singularly perturbed hybrid dynamical systems, whose stability properties are mainly characterized by the stability properties of their average hybrid dynamics. For all algorithms, structural robustness properties were established with respect \tcbb{to sufficiently small bounded disturbances.} Additionally, discretization mechanisms based on Euler and Runge-Kutta methods were also presented. Future directions will study the application of the HAES in dynamic plants, as well as adaptive and event-triggered restarting mechanisms.  In order to obtain our main results, we developed an extended averaging theorem for hybrid dynamical systems that generate an average system that renders a compact set SGPAS instead of UGAS. This result is instrumental for the analysis of hybrid extremum seeking controllers that go beyond those studied in this paper 
%
%
\section*{Acknowledgements}
\tcbb{The first author would like to thank Mihailo R. Jovanovic for insightful questions and comments that motivated part of this research.}
\bibliographystyle{plain}
\bibliography{Biblio}

\begin{thebibliography}{10}

\bibitem{KrsticBookESC}
K.~B. Ariyur and M.~Krsti\'{c}.
\newblock {\em Real-Time Optimization by Extremum-Seeking Control}.
\newblock Wiley, 2003.

\bibitem{ShivNotes}
Shivkumar Chandrasekaran.
\newblock {\em Core Matrix Analysis}.
\newblock Lecture Notes, University of California, Santa Barbara, 2010.

\bibitem{Nesterov_derivative_free}
A.~D'Aspremont.
\newblock Smooth optimization with approximate gradient.
\newblock {\em SIAM Journal of Optimization}, 19(1171-1183), 2008.

\bibitem{InexactNesterov2}
O.~Devolder, F.~Glineur, and Y.~Nesterov.
\newblock First-order methods of smooth convex optimization with inexact
  oracle.
\newblock {\em Mathematical Programming}, 146:37--75, 2014.

\bibitem{Durr13L}
H.~Durr, C.~Zeng, and C.~Ebenbauer.
\newblock Saddle point seeking in convex optimization problems.
\newblock {\em 9th IFAC Symposium on Nonlinear Control Systems}, pages
  540--545, 2013.

\bibitem{Faziyab18Siam}
M.~Faziyab, A.~Ribeiro, M.~Morari, and V.~M. Preciado.
\newblock Analysis of optimization algorithms via integral quadratic
  constraints: Nonstrongly convex problems.
\newblock {\em SIAM J. Optim.}, 28(3):2654--2689, 2018.

\bibitem{EnenbauerZeroOrder}
J.~Feiling, A.~Zeller, and C.~Ebenbauer.
\newblock Derivative-free optimization algorithms based on non-commutative
  maps.
\newblock {\em IEEE Control Systems Letters}, 2(4):743--748, 2018.

\bibitem{Newton}
A.~Ghaffari, M.~Krsti\'{c}, and D.~Ne{\v{s}}i{\'c}.
\newblock Multivariable newton-based extremum seeking.
\newblock {\em Automatica}, 48:1759--1767, 2012.

\bibitem{Goebel:12}
R.~Goebel, R.~G. Sanfelice, and A.~R. Teel.
\newblock {\em Hybrid Dynamical Systems}.
\newblock Princeton University Pressl, Princeton, NJ, USA, 2012.

\bibitem{Grushkovskaya2017}
V.~Grushkovskaya, H.~Durr, C.~Ebenbauer, and A.~Zuyev.
\newblock Extremum seeking for time-varying functions using lie bracket
  approximations.
\newblock {\em IFAC-PapersOnLine}, 50(1):5222--5528, 2017.

\bibitem{Guay:03}
M.~Guay and T.~Zhang.
\newblock Adaptive extremum seeking control of nonlinear dynamic systems with
  parametric uncertainties.
\newblock {\em Automatica}, 39:1283--1293, 2003.

\bibitem{Ronny_Kutadinata}
R.~J. Kutadinata, W.~H. Moase, and C.~Manzie.
\newblock Extremum-seeking in singularly perturbed hybrid systems.
\newblock {\em IEEE Transactions on Automatic Control}, 62(6):3014--3020, 2017.

\bibitem{NewtonESC}
C.~Labar, E.~Garone, M.~Kinnaert, and C.~Ebenbauer.
\newblock Newton-based extremum seeking: A second-order lie bracket
  approximation approach.
\newblock {\em Automatica}, 105:356--367, 2019.

\bibitem{AcceleratedMethods19}
M.~Laborde and A.~Oberman.
\newblock A {L}yapunov analysis for accelerated gradient methods: from
  deterministic to stochastic case.
\newblock In {\em In Proc. of 23rd Int. Conf. on Artificial Intelligence and
  Statistics}, volume 108, pages 602--612, 2020.

\bibitem{HeavyBollES}
S.~Michalowsky and C.~Ebenbauer.
\newblock The multidimensional n-th order heavy ball method and its application
  to extremum seeking.
\newblock {\em 53rd IEEE Conf. Decision Control}, pages 2660--2666, 2014.

\bibitem{Michalowsky15}
S.~Michalowsky and C.~Ebenbauer.
\newblock Model-based extremum seeking for a class of nonlinear systems.
\newblock {\em American Control Conference}, pages 2026--2031, 2015.

\bibitem{Mihailo19}
H.~Mohammadi, M.~Razaviyayn, and M.~R. Jovanovic.
\newblock Robustness of accelerated first-order algorithms for strongly convex
  optimization problems.
\newblock {\em IEEE Trans. Automat. Control, DOI 10.1109/TAC.2020.3008297},
  2020.

\bibitem{TAC13NesicESC}
D.~Ne{\v{s}}i{\'c}, A.~Mohammadi, and C.~Manzie.
\newblock A framework for extremum seeking control of systems with parameter
  uncertainties.
\newblock {\em {IEEE} {T}rans. {A}utom. {C}ontrol.}, 58(2):435--448, 2013.

\bibitem{DerivativesESC}
D.~Ne\u{s}i\'{c}, Y.~Tan, W.~H. Moase, and C.~Manzie.
\newblock A unifying approach to extremum seeking: Adaptive schemes based on
  estimation of derivatives.
\newblock {\em 49th IEEE Conference on Decision and Control}, pages 4625--4630,
  2010.

\bibitem{Candes_Restarting}
O'Donoghue and E.~J. Candes.
\newblock Adaptive restart for accelerated gradient schemes.
\newblock {\em Foundations of Computational Mathematics}, 15(3):715--732, 2013.

\bibitem{ThiagoKrstic}
T.~Oliveira, M.~Krsti\'{c}, and D.~Tsubakino.
\newblock Extremum seeking for static maps with delays.
\newblock {\em IEEE Trans. Autom. Control}, 62(4):1911--1926, 2017.

\bibitem{restarting_black_box}
S.~Pokutta.
\newblock Restarting algorithms: Sometimes there is free lunch.
\newblock {\em arXiv:2006.14810}, 2020.

\bibitem{BlackBoxHybridES}
J.~I. Poveda, R.~Kuttadinata, C.~Manzie, D.~Nesic, A.R. Teel, and C.~Liao.
\newblock Hybrid extremum seeking for black-box optimization in hybrid plants:
  An analytical framework.
\newblock {\em IEEE Conf. on Decision and Control}, pages 2235--2240, 2018.

\bibitem{PovedaNaliCDC19}
J.~I. Poveda and N.~Li.
\newblock Inducing uniform asymptotic stability in non-autonomous accelerated
  optimization dynamics via hybrid regularization.
\newblock {\em 58th IEEE Conference on Decision and Control}, pages 3000--3005,
  2019.

\bibitem{Poveda:16}
J.~I. Poveda and A.~R. Teel.
\newblock A framework for a class of hybrid extremum seeking controllers with
  dynamic inclusions.
\newblock {\em Automatica}, 76:113--126, 2017.

\bibitem{GuananLetters}
G.~Qu and N.~Li.
\newblock On the exponential stability of primal-dual gradient dynamics.
\newblock {\em IEEE Control Syst. Letters}, 3(1):43--48, 2019.

\bibitem{SimulatorHybridSystems}
R.~G. Sanfelice and A.~R. Teel.
\newblock Dynamical properties of hybrid systems simulators.
\newblock {\em Automatica}, 46:239--248, 2010.

\bibitem{ScheinkerBoundedLetters}
A.~Scheinker and M.~Krsti\'{c}.
\newblock Extremum seeking with bounded update rates.
\newblock {\em Systems \& Control Letters}, 63:25--31, 2014.

\bibitem{HighResolution2018}
B.~Shi, S.~Du, M.~Jordan, and W.~Su.
\newblock Understanding the acceleration phenomenon via high-resolution
  differential equations.
\newblock {\em arXiv:1810.08907}, 2018.

\bibitem{numericalanalysisbook}
A.~M. Stuart and A.~R. Humphries.
\newblock {\em Dynamical systems and numerical analysis.}
\newblock Cambridge Uniersity Press, 1996.

\bibitem{ODE_Nesterov}
W.~Su, S.~Boyd, and E.~Candes.
\newblock A differential equation for modeling {N}esterov's accelerated
  gradient method: {T}heory and insights.
\newblock {\em J. of Machine Learning Research}, 17(153):1--43, 2016.

\bibitem{Suttner2017}
R.~Suttner and S.~Dashkovskiy.
\newblock Exponential stability for extremum seeking control systems.
\newblock {\em IFAC-PapersOnLine}, 50(1):15464--15470, 2017.

\bibitem{tan06Auto}
Y.~Tan, D.~Ne\v{s}i\'{c}, and I.~M. Mareels.
\newblock On non-local stability properties of extremum seeking control.
\newblock {\em Automatica}, 42(6):889--903, 2006.

\bibitem{averagingTeel}
A.~R. Teel and D.~Ne{\v{s}}i{\'c}.
\newblock Averaging for a class of hybrid systems.
\newblock {\em Dynamics of Continuous, Discrete and Impulsive Systems},
  17:829--851, 2010.

\bibitem{HamiltonianHybrid}
A.~R. Teel, J.~I. Poveda, and J.~Le.
\newblock First-order optimization algorithms with resets and {H}amiltonian
  flows.
\newblock {\em 58th IEEE Conference on Decision and Control}, pages 5838--5843,
  2019.

\bibitem{Wang:12_Automatica}
W.~Wang, A.~Teel, and D.~Ne\u{s}i\'{c}.
\newblock Analysis for a class of singularly perturbed hybrid systems via
  averaging.
\newblock {\em Automatica}, 48(6):1057--1068, 2012.

\bibitem{Wibisono1}
A.~Wibisono, A.~C. Wilson, and M.~I. Jordan.
\newblock A variational perspective on accelerated methods in optimization.
\newblock {\em Proceedings of the National Academy of Sciences},
  113(47):E7351--E7358, 2016.

\bibitem{Ye:16b}
M.~Ye and G.~Hu.
\newblock Distributed extremum seeking for constrained network optimization and
  its application to energy consumption control in smart grid.
\newblock {\em IEEE Transactions on Control Systems Technology},
  24(6):2048--2049, 2016.

\bibitem{Jadbabaie_RK}
J.~Zhang, A.~Mokhtari, S.~Sra, and A.~Jadbabaie.
\newblock Direct runge-kutta discretization achieves acceleration.
\newblock {\em arXiv preprint arXiv:1805.00521}, 2018.

\end{thebibliography}

\appendix

\appendix

\section{Perturbed Hybrid Dynamical Systems}
\label{AppendixA}
The following proposition is a modest extension of \cite[Lem. 7.20]{Goebel:12} for the case when a nominal HDS renders a compact set SGPpAS instead of UGpAS.  
\begin{proposition}\label{robustness_proposition}
Suppose that a $\delta$-parameterized hybrid system $\mathcal{H}_{\delta}:=\{C_{\delta},D_{\delta},F_{\delta},G_{\delta}\}$ satisfies the Basic Conditions for each $\delta>0$, and that it renders a compact set $\mathcal{A}$ SGPpAS as $\delta\to0^+$ with $\beta\in\mathcal{K}\mathcal{L}$. Then, the $\rho$-inflated system $\mathcal{H}_{\delta,\rho}:=\{C_{\delta,\rho},D_{\delta,\rho},F_{\delta,\rho},G_{\delta,\rho}\}$ with data:
\begin{subequations}\label{perturbedHDS}
\begin{align}
F_{\delta,\rho}(x):&=\overline{\text{co}}~F_{\delta}((x+\rho\mathbb{B})\cap C_{\delta})+\rho\mathbb{B}\\
G_{\delta,\rho}(x):&=\{v\in\mathbb{R}^n:v\in g+\rho\mathbb{B},g\in G_{\delta}((x+\rho\mathbb{B})\cap D_{\delta})\}\\
C_{\delta,\rho}:&=\{x\in\mathbb{R}^n:(x+\rho\mathbb{B})\cap C_{\delta}\neq\emptyset\}\label{C_inflated}\\
D_{\delta,\rho}:&=\{x\in\mathbb{R}^n:(x+\rho\mathbb{B})\cap D_{\delta}\neq\emptyset\}\label{D_inflated}
\end{align}
\end{subequations}
renders the set $\mathcal{A}$ SGPpAS as $(\delta,\rho)\to0$ with $\beta\in\mathcal{K}\mathcal{L}$. 
\end{proposition}
\noindent
\textbf{Proof}: The proof is almost identical to the proof of \cite[Lem. 7.20]{Goebel:12}. Let $K\subset\mathcal{B}_{\mathcal{A}}$ and $\varepsilon>0$ be given. Since $\omega$ is continuous and grows unbounded as $x\to\text{bd}(\mathcal{B}_{\mathcal{A}})$ there exists an  $r_2>\varepsilon$ such that $K\subset\{x\in\mathcal{B}_{\mathcal{A}}:\omega(x)\leq r_2\}$. Choose $\nu=\varepsilon/4$. Then, since the system $\mathcal{H}_{\delta}$ is SGPpAS (w.r.t $\mathcal{B}_{\mathcal{A}}$) as $\delta\to0^+$ there exists $\delta^*>0$ such that for all $\delta\in(0,\delta^*)$ all solutions $x_{\delta}$ of $\mathcal{H}_{\delta}$ with $\omega(x_{\delta}(0,0))\leq r_2$  and all $(t,j)~\in\text{dom}(x_{\delta})$ the following holds:
\begin{equation}\label{firstin1}
\omega(x_{\delta}(t,j))\leq \beta(w(x_{\delta}(0,0)),t+j)+\varepsilon/4.
\end{equation}
Let $T>0$ be large enough such that  $\beta(r_2,t+j)\leq \frac{\varepsilon}{2}$, for all $t+j\geq T$.

\noindent
\textsl{Claim:} There exists a $\rho^*>0$ such that for all $\rho\in(0,\rho^*]$, all solutions $x_{\delta,\rho}$ to $\mathcal{H}_{\delta,\rho}$ with $\omega(x_{\delta,\rho}(0,0))\leq r_2$ and all $(t,j)\in\text{dom}(x_{\delta,\rho})$ the following holds: 
\begin{equation}\label{boundKL}
\omega(x_{\delta,\rho}(t,j))\leq \beta(\omega(x_{\delta,\rho}(0,0)),t+j)+\varepsilon/2.
\end{equation}
for all $t+j\leq 2T$. \QEDB
 
\noindent
By the selection of $T$ above, the claim implies that $\omega(x_{\delta,\rho}(t,j))\leq\varepsilon$, for all $2T\geq t+j\geq T$, such that $(t,j)\in\text{dom}(x_{\delta,\rho})$. We can recursively apply this argument restarting the solution and using $\varepsilon<r_2$ to get $\omega(x_{\delta,\rho}(t,j))\leq \varepsilon$ for all $(t,j)\in\text{dom}(x_{\delta,\rho})$ such that $t+j\geq T$.

To prove the claim, suppose by contradiction that there exists a sequence $\rho_i\searrow0$ and a sequence of solutions $x_{\delta,\rho_i}$ to $\mathcal{H}_{\delta,\rho_i}$ with $\omega(x_{\delta,\rho_i}(0,0,))\leq m$ and points $(t_i,j_i)\in\text{dom}(x_i)$ with $t_i+j_i\leq 2T$ such that \eqref{boundKL} does not hold: 
\begin{equation}\label{boundKL2}
\omega(x_{\delta,\rho_i}(t_i,j_i))> \beta(\omega(x_{\delta,\rho_i}(0,0)),t_i+j_i)+\varepsilon/2,
\end{equation}
Since $\omega(x_{\delta,\rho_i}(0,0,))\leq r_2$ implies that the sequence $x_{\delta,\rho_i}(0,0)$ lies in a compact subset of $\mathcal{B}_{\mathcal{A}}$, one can assume that it converges to some point in $\mathcal{B}_{\mathcal{A}}\cap (C\cup D)$. At this point, because of \eqref{firstin1}, the HDS $\mathcal{H}_{\delta}$ is pre-forward complete.  Since this implies that the sequence $x_{\delta,\rho_i}$ is locally eventually bounded \cite[Def. 5.24]{Goebel:12}, and since for each $\delta>0$ the system $\mathcal{H}_{\delta}$ satisfies the Basic Conditions, the graphical limit of the sequence $x_{\delta,\rho_i}$, denoted by $x_{\delta}$, will be a solution to $\mathcal{H}_{\delta}$. Without loss of generality we can assume that the sequence $(t_i,j_i)$ also converges to some $(t,j)\in\text{dom}(x_{\delta})$. Using continuity of $\omega$ and $\beta$ and taking the limit as $i\to\infty$ at both sides of \eqref{boundKL2} we obtain $\omega(x_{\delta}(t,j))> \beta(\omega(x_{\delta}(0,0)),t+j)+\frac{\varepsilon}{2}$, which violates \eqref{firstin1} at the time $(t,j)\in\text{dom}(x_{\delta})$. This is a contradiction. \null \hfill \null $\blacksquare$

\begin{lemma}\label{lemma_optimal_set}
Let $\phi:\mathbb{R}^n\to\mathbb{R}$ satisfy Assumption \ref{assumption1} with item (a), and let $\mathcal{A}_{\phi}=\{z^*\}$ be the unique minimizer of $\phi$. Consider the set
\begin{equation*}
\mathcal{O}:=\Big\{x_1\in\mathbb{R}^n:(x_1-z^*)^\top\nabla \phi(x_1) -(\phi(x_1)-\phi^*)=0\Big\}.
\end{equation*}
Then,  $\mathcal{O}\subset \mathcal{A}_{\phi}$.
\end{lemma}
\noindent
\textbf{Proof}: Let $x_1$ be such that $(x_1-z^*)^\top\nabla \phi(x_1) -(\phi(x_1)-\phi^*)=0$. Suppose that $x_1\notin\mathcal{A}_{\phi}$. Let $\alpha_1:=\phi(x_1)$ and define the set $\Omega_{\alpha_1}:=\{x\in\mathbb{R}^n: \phi(x)\leq \alpha_1\}$. Since $z^*$ is optimal, we have that $\phi(z^*)=\phi^*\leq \phi(x_1)$ and therefore $\mathcal{A}_{\phi}\subset  \Omega_{\alpha_1}$. Since $\phi\in\mathcal{C}^2$ we have that $\nabla \phi$ is locally Lipschitz, and since $\Omega_{\alpha_1}$ is compact due to the radial unboundedness of $\phi$ there exists $L_{\alpha_1}>0$ such that $|\nabla \phi(x'_1)-\nabla \phi(x''_1)|\leq L_{\alpha_1} |x'_1-x''_1|$ for all $(x'_1,x''_1)\in \Omega_{\alpha_1}$. By the convexity and Lipschitz properties in $\Omega_{\alpha_1}$:
\begin{equation*}
\phi(z^*)-\phi(x_1)-\nabla \phi(x_1)^{\top}(z^*-x_1)\geq \frac{1}{2L_{\alpha_1}}|\nabla \phi(x_1)|^2\geq0,
\end{equation*}
but since by assumption the left hand side of the inequality is zero, we must have that $\nabla \phi(x_1)=0$, which is a contradiction given that $x_1\notin \mathcal{A}_{\phi}$ and $\phi$ is convex. \null \hfill \null $\blacksquare$
\section{Solutions to Hybrid Dynamical Systems}
\label{AppendixB}
A hybrid dynamical system is modeled by the equation:
\begin{subequations}\label{HDS_appendix}
\begin{align}
x&\in C,~~~~~~~~\dot{x}= F(x),\\
x&\in D,~~~~~x^+\in G(x),
\end{align}
\end{subequations}
where $x\in\mathbb{R}^n$ is the state, and the mappings $F:\mathbb{R}^n\to\mathbb{R}^n$ and $G:\mathbb{R}^n\rightrightarrows\mathbb{R}^n$ and the sets $C\subset\mathbb{R}^n$ and $D\subset\mathbb{R}^n$ satisfy the Basic Conditions \cite[Assumption 6.5]{Goebel:12}. A HDS that satisfies the Basic Conditions is said to be \emph{well-posed} \cite[Thm. 6. 30]{Goebel:12}, which permits the use of graphical convergence tools to establish sequential compactness results for the solutions of \eqref{HDS_appendix} (e.g., the graphical limit of a sequence of solutions is also a solution) see \cite[Sec. 6.2-6.4]{Goebel:12}. When the jump map is single-valued, e.g., as in \eqref{HDS}, one can define the set-valued map $G$ as $G(x)=G_0$ when $x\in D$, and $G(x)=\emptyset$ when $x\notin D$, where $x\mapsto G_0(x)$ is the original single-valued function. Solutions of \eqref{HDS_appendix} are defined on \emph{hybrid time domains} \cite[Ch. 2]{Goebel:12}. A set $E\subset\mathbb{R}_{\geq0}\times\mathbb{Z}_{\geq0}$ is called a \textsl{compact} hybrid time domain if $E=\cup_{j=0}^{J-1}([t_j,t_{j+1}],j)$ for some finite sequence of times $0=t_0\leq t_1\ldots\leq t_{J}$. The set $E$ is a hybrid time domain if for all $(T,J)\in E$, the set $E\cap([0,T]\times\{0,\ldots,J\})$ is a compact hybrid time domain. 
\begin{definition}\cite[pp. 124]{Goebel:12}
A function $x:\text{dom}(x)\mapsto\mathbb{R}^n$ is a \emph{hybrid arc} if $\text{dom}(x)$ is a hybrid time domain and $t\mapsto x(t,j)$ is locally absolutely continuous for each $j$ such that the interval $I_j:=\{t:(t,j)\in \text{dom}(x)\}$ has nonempty interior. A hybrid arc $x$ is a \emph{solution} to a well-posed HDS \eqref{HDS_appendix} if $x(0,0)\in C\cup D$, and the following two conditions hold: (1): For each $j\in\mathbb{Z}_{\geq0}$ such that $I_j$ has nonempty interior: $x(t,j)\in C$ and $\dot{x}(t,j)\in F(x(t,j))$ for almost all $t\in I_j$; (2): For each $(t,j)\in\text{dom}(x)$ such that $(t,j+1)\in \text{dom}(x)$: $x(t,j)\in D$, and $x(t,j+1)\in G(x(t,j))$.
\end{definition}
The following definition and lemma characterize maximal and complete solutions in HDS.
\begin{definition}\cite[Def. 2.5 $\&$ 2.7]{Goebel:12}
A hybrid solution $x$ is said to be: a) nontrivial if $\text{dom}(x)$ contains at least two points; b) forward pre-complete if its domain is compact or unbounded; c) complete if its hybrid time domain is unbounded; d) maximal if there does not exist another solution $\psi$ to $\mathcal{H}$ such that $\text{dom}(x)$ is a proper subset of $\text{dom}(\psi)$, and $x(t,j)=\psi(t,j)$ for all $(t,j)\in\text{dom}(x)$. 
\end{definition}
\begin{lemma}\cite[Prop. 6.10]{Goebel:12}\label{existence_completness}
Let $\mathcal{H}=(C,F,D,G)$ satisfy the Basic Conditions \cite[Assumption 6.5]{Goebel:12}. Let $T_X(y)$ be the tangent cone to a set $X$ at a point $y$. Let $x_0\in C\cup D$. If either (a) $x_0\in D$; or (b) there exists a neighborhood $U$ of $x_0$ such that for every $y\in U\cap C$ we have that $F(y)\cap T_C(y)\neq\emptyset$, then there exists a nontrivial solution $x$ to $\mathcal{H}$ with $x(0,0)=x_0$. Moreover, if item (b) holds for every $x_0\in C\backslash D$, then there exists a nontrivial solution to $\mathcal{H}$ from every initial point in $C\cup D$, and every maximal solution $x$ satisfies exactly one of the following conditions: a) $x$ is complete; b) $x$ has a finite escape time; c) $x(T,J)\notin C\cup D$, where $(T,J)=\sup~\text{dom}(x)$, i.e., the solution $x$ stops. Furthermore, if $G(D)\subset C\cup D$, then (c) above does not occur.
\end{lemma}

\section{Proof of Lemma \ref{lemma_dither}}
\label{sec_lemma_dither}
Let $\Psi(t):=\left[\Psi_1(t)^\top,\Psi_3(t)^\top,\Psi_5(t)^\top,\ldots,\Psi_{2n-1}(t)^\top\right]^\top\in\mathbb{R}^{2n}$ and $\mu_{0}:=[\mu_{0,1}^\top,\mu_{0,3}^\top,\mu_{0,5}^\top,\ldots,\mu_{0,n-1}^\top]^\top\in\mathbb{R}^{2n}$, where the vectors $\Psi_1(t)$ and $\mu_{0,i}(t)$ are defined as in \eqref{sinusoid2} with $1/\varepsilon=1$. With these definitions in hand, we can write $\tilde{\mu}(t)=\text{diag}(\mu_0)^\top\Psi(t)$ which implies \cite[Exercise 260]{ShivNotes} that
\begin{equation*}
\int \tilde{\mu}(t)\tilde{\mu}(t)^\top dt=\text{diag}(\mu_0)^\top\left(\int\Psi(t)\Psi(t)^\top dt\right) \text{diag}(\mu_0),
\end{equation*}
where $\text{diag}(\mu_0)$ is a $(2n\times n)$ block diagonal matrix with diagonal blocks given by $\mu_{0,i}$, $i\in\{1,3,5,\ldots,n-1\}$. Thus, it suffices to show the existence of a $\tilde{T}>0$ such that 
\begin{equation}\label{important_property}
\frac{1}{k\tilde{T}}\int_{0}^{k\tilde{T}}\Psi(t)\Psi(t)^\top dt=\frac{1}{2}I_{2n},~~~\int_{0}^{k\tilde{T}}\Psi(t)=\mathbf{0}_{2n},
\end{equation}
since this would imply that $\int_{0}^{k\tilde{T}} \tilde{\mu}(t)dt=\textbf{0}_n$ and that $\int_{0}^{k\tilde{T}} \tilde{\mu}(t)\tilde{\mu}(t)^\top dt=0.5\text{diag}(\mu_0)^\top\text{diag}(\mu_0)=0.5I_n$, where the last equality follows by the fact that $\mu_{0,i}^\top\mu_{0,i}=1$ for all $i\in\{1,3,\ldots,n\}$ since $\mu(0)\in\mathbb{T}^n$. To show \eqref{important_property} we show the existence of a $\tilde{T}>0$ such that
\begin{equation}\label{proof_property}
\frac{1}{k\tilde{T}}\int_{0}^{k\tilde{T}}\Psi_i(t)\Psi_j(t)^\top dt=c_{ij} I_{2},~~\int_{0}^{k\tilde{T}}\Psi_i(t)=\mathbf{0}_2,
\end{equation}
for all $i,j\in\{1,3,5,\ldots,n\}$ and all $k\in\mathbb{Z}_{>0}$, where $c_{ij}=0.5$ for all $i=j$, and $c_{ij}=0$ for all $i\neq j$. Indeed, by Assumption \ref{assumption_freq}, the parameters $\kappa_{\ell}$ can be written as $\kappa_{\ell}=\kappa_{\ell}^n/\kappa^d_{\ell}$, for all $\ell\in\{1,2,3,\ldots,n\}$, where $\kappa^n_{\ell}$ and $\kappa_{\ell}^d$ are positive integers. Let $T_{\ell}:=1/\kappa_{\ell}$ and $\kappa^n:=\Pi_{j=1}^n\kappa^n_{j}\in\mathbb{Z}_{>0}$. Then, $T_{\ell} \kappa^n=\kappa^d_{\ell}\Pi_{j=1,j\neq\ell}^n\kappa^n_{j}\in\mathbb{Z}_{>0}$. Define $\tilde{T}_{\ell}:=T_{\ell}\kappa^n$ and let $\tilde{T}=\text{LCM}\{\tilde{T}_1,\tilde{T}_2,\ldots,\tilde{T}_{n}\}$, where $\text{LCM}$ stands for least common multiplier. Such $\tilde{T}$ is a well defined positive integer and it is unique. Then, by definition of the $\text{LCM}$, for each $\ell$ there exists a $n_{\ell}\in\mathbb{Z}_{>0}$ such that 
\begin{equation}\label{LCM}
\tilde{T}=n_{\ell}\tilde{T}_{\ell}=\tilde{n}_{\ell} T_{\ell},~~~\tilde{n}_{\ell}:=n_{\ell}\kappa^n\in\mathbb{Z}_{>0}.
\end{equation}
%
%
%
Using $\ell=(i+1)/2$ and the definition of $\Psi_i(t)$, we obtain
\begin{equation*}
\int_{0}^{k\tilde{T}}\Psi_i(t)dt=\left[\begin{array}{c}
\int_{0}^{k\tilde{T}}\cos\left(\frac{2\pi}{T_{\ell}}t\right)dt\\
\int_{0}^{k\tilde{T}}\sin\left(\frac{2\pi}{T_{\ell}}t\right)dt
\end{array}\right]=\left[\begin{array}{c}
\sin\left(\frac{2\pi}{T_{\ell}}t\right)\Big|_{0}^{k\tilde{n}_{\ell}T_{\ell}}\\
-\cos\left(\frac{2\pi}{T_{\ell}}t\right)\Big|_{0}^{k\tilde{n}_{\ell}T_{\ell}}
\end{array}\right],
\end{equation*}
which is equal to $\mathbf{0}_2$ for all $i\in\{1,3,5,\ldots,2n-1\}$, $k\in\mathbb{Z}_{\geq0}$, as in \eqref{proof_property}. Also, $\int_{0}^{k\tilde{T}}\Psi_i(t)\Psi_j(t)^\top dt$ is
\begin{equation}\label{matrix_dithersignal}
=\left[\begin{array}{cc}
\int_{0}^{k\tilde{T}}\cos\left(\frac{2\pi}{T_{\ell}}t\right)\cos\left(\frac{2\pi}{T_{s}}t\right)dt &~ \int_{0}^{\tilde{kT}}\cos\left(\frac{2\pi}{T_{\ell}}t\right)\sin\left(\frac{2\pi}{T_{s}}t\right)dt\\
\int_{0}^{k\tilde{T}}\sin\left(\frac{2\pi}{T_{\ell}}t\right)\cos\left(\frac{2\pi}{T_{s}}t\right)dt &~ \int_{0}^{k\tilde{T}}\sin\left(\frac{2\pi}{T_{\ell}}t\right)\sin\left(\frac{2\pi}{T_{s}}t\right)dt
\end{array}\right],
\end{equation}
where $s=(j+1)/2$. When $i=j$ we have that $\ell=s$ and the diagonal terms satisfy
\begin{align*}
\int_{0}^{k\tilde{T}}\cos\left(\frac{2\pi}{T_{\ell}}t\right)^2dt&=\frac{1}{2}\left(t+\frac{\sin(\frac{4\pi}{T_{\ell}}t)T_{\ell}}{4\pi}\right)\Bigg|_{0}^{k\tilde{n}_{\ell}T_{\ell}} =\frac{k\tilde{n}_{\ell}T_{\ell}}{2}\\
\int_{0}^{k\tilde{T}}\sin\left(\frac{2\pi}{T_{\ell}}t\right)^2dt&=\frac{1}{2}\left(t-\frac{\sin(\frac{4\pi}{T_{\ell}}t)T_{\ell}}{4\pi}\right)\Bigg|_{0}^{k\tilde{n}_{\ell}T_{\ell}} =\frac{k\tilde{n}_{\ell}T_{\ell}}{2},
\end{align*}
while the off-diagonal terms are given by
\begin{align*}
\int_{0}^{k\tilde{T}}\cos\left(\frac{2\pi}{T_{\ell}}t\right)\sin\left(\frac{2\pi}{T_{\ell}}t\right)dt&=\frac{\sin\left(\frac{2\pi}{T_{\ell}}t\right)^2T_{\ell}}{4\pi}\Bigg|_{0}^{k\tilde{n}_{\ell}T_{\ell}}=0.
\end{align*}
Thus, when $i=j$ we have that $c_{ij}=0.5$ in \eqref{proof_property}. On the other hand, when $i\neq j$, we have that $\ell \neq s$ and the diagonal terms of \eqref{matrix_dithersignal} become
\begin{align*}
\int_{0}^{k\tilde{T}}\cos\left(\frac{2\pi}{T_{\ell}}t\right)\cos\left(\frac{2\pi}{T_{s}}t\right)dt&=\frac{\sin(2\pi tT^+_{\ell,s})}{4\pi(T^+_{\ell,s})}+\frac{\sin(2\pi tT^-_{\ell,s})}{4\pi(T^-_{\ell,k})}\Bigg|^{k\tilde{T}}_{0}\\
\int_{0}^{k\tilde{T}}\sin\left(\frac{2\pi}{T_{\ell}}t\right)\sin\left(\frac{2\pi}{T_{s}}t\right)dt&=\frac{\sin(2\pi tT^-_{\ell,s})}{4\pi(T^-_{\ell,s})}-\frac{\sin(2\pi tT^+_{\ell,s})}{4\pi(T^+_{\ell,s})}\Bigg|^{k\tilde{T}}_{0}
\end{align*}
where $T^+_{\ell,s}=\frac{1}{T_\ell}+\frac{1}{T_s}$ and $T^-_{\ell,s}=\frac{1}{T_\ell}-\frac{1}{T_s}$. Using \eqref{LCM}, when $t=k\tilde{T}$ we get $tkT^+_{\ell,s}=k(\tilde{n}_{\ell}+\tilde{n}_{s})\in\mathbb{Z}_{>0}$. Similarly, $tkT^-_{\ell,s}=k(\tilde{n}_{\ell}-\tilde{n}_{s})\in\mathbb{Z}_{\neq0}$. This implies that both integrals are zero for any $k\in\mathbb{Z}_{>0}$. Finally, when $s\neq\ell$ the off-diagonal terms satisfy
\begin{align*}
\int_{0}^{k\tilde{T}}\sin\left(\frac{2\pi}{T_{\ell}}t\right)\cos\left(\frac{2\pi}{T_{s}}t\right)dt&=-\frac{\cos(2\pi tT^-_{\ell,s})}{4\pi(T^-_{\ell,s})}-\frac{\cos(2\pi tT^+_{\ell,s})}{4\pi(T^+_{\ell,s})}\Bigg|^{k\tilde{T}}_{0}\\
\int_{0}^{k\tilde{T}}\cos\left(\frac{2\pi}{T_{\ell}}t\right)\sin\left(\frac{2\pi}{T_{s}}t\right)dt&=\frac{\cos(2\pi tT^-_{\ell,s})}{4\pi(T^-_{\ell,s})}-\frac{\cos(2\pi tT^+_{\ell,s})}{4\pi(T^+_{\ell,s})}\Bigg|^{k\tilde{T}}_{0}
\end{align*}
which are also zero by the definition of $T^+_{\ell,s}$, $T^-_{\ell,s}$, and \eqref{LCM}. This establishes that $c_{ij}=0$ in \eqref{proof_property} whenever $i\neq j$. \null \hfill \null $\blacksquare$

\end{document}